\documentclass[12pt,oneside]{amsart}
\usepackage[a4paper,left=2cm,top=3cm,right=2cm,bottom=2cm]{geometry}
\usepackage{amsmath,amssymb,amsthm,mathtools}
\usepackage{enumitem}
\usepackage{mathrsfs}
\usepackage[dvipsnames]{xcolor}
\usepackage[
  colorlinks,
  citecolor=blue,
  linkcolor=blue,
  urlcolor=blue
]{hyperref}
\usepackage[utf8]{inputenc}
\usepackage[T1]{fontenc}
\expandafter\providecommand\csname subjclassname@2020\endcsname{\textup{2020} Mathematics Subject Classification}
\theoremstyle{plain}
\newtheorem{theorem}{Theorem}[section]
\newtheorem{proposition}[theorem]{Proposition}
\newtheorem{lemma}[theorem]{Lemma}
\newtheorem{corollary}[theorem]{Corollary}
\newtheorem*{theoremA}{Theorem A}
\newtheorem*{theoremB}{Theorem B}
\newtheorem*{theoremC}{Theorem C}
\newtheorem*{theoremE}{Theorem E}
\newtheoremstyle{introcorollary}{6pt}{6pt}{\itshape}{}%
  {\bfseries}{.}{.5em}{}
\theoremstyle{introcorollary}
\newtheorem*{corollaryD}{Corollary D}
\theoremstyle{definition}
\newtheorem{definition}[theorem]{Definition}
\newtheorem{example}[theorem]{Example}
\theoremstyle{remark}
\newtheorem{remark}[theorem]{Remark}
\numberwithin{equation}{section}
\newcommand{\F}{\mathcal F}
\newcommand{\Fl}{\mathcal F^\ell}
\newcommand{\G}{\mathcal G}
\newcommand{\Gl}{\mathcal G^\ell}
\newcommand{\R}{\mathbb R}
\newcommand{\C}{\mathbb C}
\newcommand{\Z}{\mathbb Z}
\newcommand{\Oct}{\mathbb O}
\newcommand{\kk}{\mathfrak k}
\newcommand{\so}{\mathfrak{so}}
\newcommand{\Hol}{\mathrm{Hol}}
\newcommand{\Fr}{\mathrm{O}}

\newcommand{\rra}{\rightrightarrows}
\newcommand{\metric}{\ensuremath{\mathrm{g}}}

\title[Reeb stability for singular Riemannian foliations]
{Reeb stability for singular Riemannian foliations}

\author{Mateus de Melo}
\address{{Departamento de Matem\'atica,
Universidade Federal do Esp\'irito Santo, Vit\'oria, Brazil}}
\email{mateus.melo@ufes.br}

\author{Laura Ribeiro dos Santos}
\address{Instituto de Matem\'atica e Estat\'istica,
Universidade de S\~ao Paulo, Rua do Mat\~ao 1010,
05508-090 S\~ao Paulo, Brazil}
\email{laura.ribeiro.santos@usp.br}

\author{Ivan Struchiner}
\address{Instituto de Matem\'atica e Estat\'istica,
Universidade de S\~ao Paulo, Rua do Mat\~ao 1010,
05508-090 S\~ao Paulo, Brazil}
\email{ivanstru@ime.usp.br}


\begin{document}

\begin{abstract}
 We prove a global Reeb stability theorem for singular Riemannian foliations on compact connected manifolds using the semi-local model developed in \cite{AIMS2022}. We show that all leaves are compact if and only if the foliation is infinitesimally closed, and there exists a compact leaf with a linear holonomy groupoid which is proper at that leaf. These conditions are also equivalent to the existence of a smooth map to a Euclidean space whose fibers are exactly the leaves. The proof combines a local criterion for strong stability with compact transverse leaves, Molino's structural theory, and polynomial descriptions of local leaf spaces. Around a closed leaf, we show that the linear holonomy groupoid realizes the original foliation precisely when the infinitesimal foliation is homogeneous. When the infinitesimal foliation is closed, properness of this groupoid yields an associated-bundle model whose compact structure group is its isotropy group. We also characterize the existence, near a compact leaf, of a Hausdorff model Lie groupoid which is proper at that leaf. As an application, we obtain stability results for maps whose connected fiber components form singular Riemannian foliations. We compare linear holonomy with Fernandes's reduced Lie algebroid holonomy \cite{Fernandes2002} and with formal outer holonomy. Finally, we prove that, for every infinitesimally closed singular Riemannian foliation, the compactly supported tangent vector fields form a locally finitely generated involutive module inducing the original leaf partition.
\end{abstract}

\maketitle

\section{Introduction}

Reeb's global stability theorem shows how the topology of one compact leaf can determine an entire foliation: for a transversely orientable
codimension-one foliation of a compact connected manifold, the existence of a compact leaf with finite fundamental group forces all leaves to be compact \cite{Reeb1952} (see also \cite[Section~2.5]{MoerdijkMrcun2003}). Its local counterpart gives a saturated neighbourhood of compact leaves around a compact leaf with finite holonomy \cite[Section~2.3]{MoerdijkMrcun2003}. We establish a global stability criterion for singular Riemannian foliations, in arbitrary codimension, by combining a principal bundle model around closed leaves, obtained from a semi-local groupoid, with Molino's structural theory.

Let $\F$ be a singular Riemannian foliation. We call $\F$ \emph{closed} if all its leaves are closed in the ambient manifold. Around a closed leaf $L$, fix $p\in L$ and a tube $U$ carrying the semi-local model, and write $\F_p$ for the infinitesimal foliation on $V=\nu_pL$.
We write $B_\varepsilon(V)$ for the open ball of radius $\varepsilon$ centered at the origin in $V$. The semi-local construction of \cite[Theorems~1.4 and~4.1]{AIMS2022} associates to this tube a Lie groupoid $\Gl$, the \emph{linear holonomy groupoid}. Its orbits
are the leaves of the linearized foliation $\Fl\subset\F$. Our global theorem relates properness of this groupoid at a single compact leaf to compactness of every leaf and to a smooth presentation of the foliation.

A leaf is \emph{stable} if it admits arbitrarily small saturated open neighbourhoods. A compact leaf $L$ is \emph{strongly stable} if it
admits the following local model (Definition~\ref{def:reeb-stability}).
A saturated neighbourhood of $L$ is foliated diffeomorphic to an
associated bundle $P\times_H\Omega$, with $L$ as the zero section.
Here $P\to L$ is a principal bundle with compact structure group $H$,
and $\Omega$ is a foliated neighbourhood of the origin in a linear
$H$-representation. The foliation on the bundle descends from the
product foliation on $P\times\Omega$. We say that the model has
\emph{compact transverse leaves} if the leaves in $\Omega$ are compact.

\begin{theoremA}[Global Reeb stability (Theorem~\ref{thm:global})]
Let $(M,\F)$ be a singular Riemannian foliation on a compact connected
manifold. The following are equivalent:
\begin{enumerate}[label=(\roman*)]
\item There is a compact leaf $L$ such that $\Gl$ is proper at $L$
and $\F_p$ is closed for some $p\in L$.
\item Some leaf of $\F$ admits a strong stability model with compact
transverse leaves.
\item Every leaf of $\F$ is compact.
\item The leaves are exactly the fibers of a smooth map $M\to\R^N$.
\end{enumerate}
\end{theoremA}

Under these equivalent conditions, every leaf is strongly stable.
In particular, a compact leaf with finite fundamental group and closed
infinitesimal foliation forces all leaves to be compact
(Corollary~\ref{cor:global-finite}). For regular Riemannian foliations,
the infinitesimal foliation is by points, giving a finite-fundamental-group
criterion analogous to classical global Reeb stability. In the singular
case the closedness hypothesis is essential: an irrational isometric
flow on $S^4$ has point leaves but also nonclosed leaves
(Example~\ref{ex:irrational-sphere}).

The local stability theorem shows that (i) and (ii) are equivalent.
Both amount to the existence of a compact leaf with a saturated tube
consisting of compact leaves. Molino's locally constant structural
sheaf then implies that every leaf of $\F$ is closed
(Lemma~\ref{lem:propagation}). Since $M$ is compact, this gives (iii).
Finally, we use polynomial presentations of the local leaf spaces and
basic cutoff functions to construct the smooth map in (iv).

The passage from (i) to (ii) has a classical counterpart. For a
regular foliation, a compact leaf $L$ with finite holonomy group $H$
has a saturated neighbourhood modeled on
$\widetilde L\times_H B_\varepsilon(V)$, where $\widetilde L$ is
the holonomy cover and the ball in $V=\nu_pL$ is taken with
respect to an $H$-invariant inner product. This is the linear
model of the holonomy groupoid along $L$. For linearization of Lie
groupoids, see
\cite{Weinstein2002,Zung2006,CrainicStruchiner2013,delHoyoFernandes2018}.

For a regular Riemannian foliation, the linear holonomy groupoid
along a closed leaf is precisely the groupoid underlying this model:
$\Gl=\Hol(\F)|_L\ltimes\nu L$ \cite[Example~4.5]{AIMS2022}.
For a singular Riemannian foliation, the orbits of $\Gl$ are the
leaves of $\Fl$, obtained by linearizing vector fields tangent to
$\F$ along $L$ \cite[Section~2]{AlexandrinoRadeschi2017},
\cite[Sections~2.1.1--2.1.2 and Theorem~4.1]{AIMS2022}.
A leaf of $\F$ can contain several such orbits.

To prove Theorem~A we need to know when the leaves of $\F$ near a
compact leaf are compact. We first ask when they are closed in a
tube around $L$, allowing $L$ itself to be any closed leaf.

\begin{theoremB}[Theorem~\ref{thm:closed}]
Let $L$ be a closed leaf, let $p\in L$, and let $U$ be a tube around
$L$ carrying the semi-local model. Write $V=\nu_pL$.
\begin{enumerate}[label=(\roman*)]
\item \emph{Closedness.}
The linearized foliation $\Fl|_U$ is closed if and only if $\Gl|_U$
is proper. The original foliation $\F|_U$ is closed if and only if
both $\Fl|_U$ and $\F_p$ are closed.
\item \emph{Local model.}
If $\F|_U$ is closed, there exist a compact Lie group $H$, a principal
$H$-bundle $P\to L$, and a foliated diffeomorphism
\[
(U,\F|_U)\cong P\times_H(B_\varepsilon(V),\F_p).
\]
If $L$ is compact, this model makes $L$ strongly stable and hence
stable.
\end{enumerate}
\end{theoremB}

Here $H$ is the isotropy group of $\Gl$ at $0_p$, and $P$ is the
source fiber over $p$ of the restriction of $\Gl$ to $L$. The action
of $H$ on $V$ is faithful. Its identity component acts as $K_p^0$,
the largest connected group of linear isometries preserving each leaf
of $\F_p$ (Lemma~\ref{lem:structure}). Properness of $\Gl|_U$ is
equivalent to compactness of $H$, or equivalently to closedness of
$K_p^0$ in $O(V)$ and finiteness of $\Gamma=\pi_0(H)$
(see~\eqref{eq:proper-isotropy}). When $\F_p$ is closed, $K_p^0$ is
compact (Lemma~\ref{lem:closed-K}). Since $\Gamma$ is a quotient of
$\pi_1(L,p)$, finite fundamental group then suffices for closedness
of $\F|_U$ (Corollary~\ref{cor:finite-components}).

The description of $H$ as an isotropy group also clarifies how
Theorem~B relates to the slice theorem of Mendes and Radeschi
\cite[Theorem~A]{MendesRadeschi2019}. Both give an associated-bundle
model for the tube. The structure group in the slice theorem consists
of all orthogonal transformations of $V$ preserving each leaf of the
disconnected slice foliation, whose leaves are the intersections of
the leaves of $\F$ with the normal slice. The action of $H$ on $V$
identifies it with a possibly proper subgroup of this group.

Suppose now that $L$ is compact. Each leaf in the tube lies at a
fixed distance from $L$, hence in a compact sphere bundle over $L$
or in $L$ itself. Closedness in the tube is therefore equivalent to
compactness. In this case the model in Theorem~B has compact
transverse leaves, giving the local criterion used in Theorem~A.

\begin{theoremC}[Theorem~\ref{thm:reeb}]
Let $L$ be a compact leaf and $p\in L$. After shrinking the saturated
tube $U$ around $L$, the following are equivalent:
\begin{enumerate}[label=(\roman*)]
\item $L$ admits a strong stability model with compact transverse leaves.
\item Every leaf of $\F|_U$ is compact.
\item $\Gl$ is proper at $L$ and $\F_p$ is closed.
\end{enumerate}
\end{theoremC}

Under these conditions, $L$ is stable. The principal bundle $P\to L$
of Theorem~B has compact total space, and $\Gl|_U$ is Hausdorff and
$s$-proper. In the regular case, $\F_p$ is the foliation by points,
$H$ is the holonomy group, and $P$ is the holonomy cover. Thus
Theorem~C recovers classical local Reeb stability
(Example~\ref{rem:classical}; see also
\cite[Example~1.15]{CrainicStruchiner2013}). Together with
Corollary~\ref{cor:finite-components}, it gives the
finite-fundamental-group consequence of Theorem~A. Related local
stability results for compact transversally stable leaves were
obtained by Dazord \cite{Dazord1985}.

The principal bundle model also yields finite covering and leaf-bundle
descriptions (Corollaries~\ref{cor:finite-cover}--\ref{cor:leaf-bundle}),
as well as polynomial presentations of the local leaf spaces and smooth
factorization of basic functions (Corollary~\ref{cor:polynomial}).
The latter uses the slice and invariant-function methods of
Mendes--Radeschi \cite{MendesRadeschi2019} and Lytchak--Radeschi
\cite{LytchakRadeschi2018} and supplies the local maps needed in
Theorem~A(iv).

In the bundle model of Theorem~B, the leaves of $\F|_U$ are the
images of $P\times S$, where $S$ ranges over the leaves of $\F_p$
in $B_\varepsilon(V)$. The orbits of $\Gl$ are the leaves of $\Fl|_U$.
When do these two foliations coincide?

On a normal fiber, the infinitesimal leaves of $\Fl$ are the
$K_p^0$-orbits \cite[Proposition~2.10]{AlexandrinoRadeschi2017}.
For an orbit-like foliation, the infinitesimal leaves of $\F$ are
compact and homogeneous \cite{Molino1994}, so the two descriptions
agree. More generally, if $B$ is a closed saturated submanifold
contained in a stratum, then $\F|_U=\Fl|_U$ on a model tube around
$B$ precisely when each reduced infinitesimal foliation $\F_b$,
$b\in B$, is homogeneous (Proposition~\ref{prop:linearizable}).
Along a closed leaf $L$, homogeneity at one point suffices. It also
characterizes the existence of a foliated diffeomorphism between
$\F$ and $\Fl$ near $L$ carrying $L$ onto itself
(Corollary~\ref{cor:diffeomorphic}). The connected groups realizing
homogeneity need not be closed. Thus homogeneity gives a Lie groupoid
model even without properness. In this case Theorem~B identifies
properness with closedness of the original leaves
(Corollary~\ref{cor:homogeneous-proper}).

When the leaves in the tube are closed, Theorem~B gives an
associated-bundle model for $\F$ whose compact structure group is the
isotropy of $\Gl$. One might therefore ask whether the action of this
group already determines the transverse leaves as its orbits. In the
octonionic Hopf foliation on $\R^{16}$, all leaves are compact, but
the isotropy at the origin is trivial: every linear vector field
tangent to the leaves vanishes \cite[Lemma~4.19]{NahariStrobl2024}.
Thus $\Gl$ is the unit groupoid, while the leaves away from the origin
are $7$-spheres (Example~\ref{ex:octonionic-linearization}).

This leaves open the possibility that a different Lie groupoid might
realize the same foliation. For the Hopf example such groupoids do
exist \cite[Theorem~3.3]{NahariStrobl2024}. The following criterion
explains when one can choose a Hausdorff Lie groupoid which is proper
at a given leaf.

\begin{corollaryD}[Corollary~\ref{cor:proper-model}]
Let $L$ be a compact leaf and $p\in L$. Near $L$, the foliation $\F$
is the orbit foliation of a Hausdorff Lie groupoid which is proper
at $L$ if and only if $L$ admits a strong stability model with compact
transverse leaves and $\F_p$ is homogeneous. In this case $\F=\Fl$
near $L$, and $\Gl$ provides such a model.
\end{corollaryD}

For a continuous map to a Hausdorff space whose connected fiber
components form a singular Riemannian foliation, closedness is
automatic. Our criterion gives strong stability at each compact fiber
component. At such a component, homogeneity of the infinitesimal
foliation characterizes the existence of a Hausdorff Lie groupoid
model which is proper at that component (Corollary~\ref{cor:fibers}).
For a general closed saturated centre, we obtain a properness
criterion and a sufficient condition for the existence of a proper
model (Proposition~\ref{prop:proper-B} and
Corollary~\ref{cor:proper-B-model}).

We now ask what the bundle model tells us about the vector fields
tangent to its leaves. Androulidakis--Skandalis associate a
\emph{topological} holonomy groupoid to every locally finitely
generated involutive module of vector fields
\cite{AndroulidakisSkandalis2009}. Androulidakis--Zambon relate
linearizability of such a module around a leaf to compactness of
holonomy isotropy and to proper Lie groupoid models
\cite[Definition~4.9 and Proposition~4.12]{AndroulidakisZambon2014}.
To describe the module in our bundle model, including when the
infinitesimal foliation is inhomogeneous, we first establish local
finite generation for the module of all leaf-tangent fields.
Recall that $\F$ is \emph{infinitesimally closed} if every $\F_p$
has closed leaves. The following theorem extends Nahari--Strobl's
result for the octonionic Hopf foliation
\cite[Proposition~4.18]{NahariStrobl2024}.

\begin{theoremE}[Theorem~\ref{thm:maximal-module}]
Let $\F$ be an infinitesimally closed singular Riemannian foliation.
Then
\[
\mathfrak X_c(\F)
=
\{X\in\mathfrak X_c(M):X(x)\in T_xL_x
\text{ for every }x\in M\}
\]
is a locally finitely generated involutive $C^\infty(M)$-module
whose leaf partition is $\F$. It is maximal among modules realizing
this partition and is natural under foliated diffeomorphisms.
\end{theoremE}

Under the equivalent conditions of Theorem~C, the bundle model gives
an explicit description of $\mathfrak X_c(\F|_U)$
(Corollary~\ref{cor:module-reeb}). In a bundle chart
$U_0\times B_\varepsilon(V)$, choose coordinates $u_1,\ldots,u_d$
on $U_0\subset L$ and local generators $Y_1,\ldots,Y_r$ for the
module of smooth fields tangent to $\F_p$. The leaf-tangent fields
on the model are locally exactly the fields of the form
\[
\sum_{i=1}^d a_i(u,v)\frac{\partial}{\partial u_i}
+\sum_{j=1}^r b_j(u,v)Y_j(v),
\]
with smooth coefficients $a_i$ and $b_j$. Under the foliated
diffeomorphism of Theorem~B, the module $\mathfrak X_c(\F|_U)$
corresponds to the compactly supported fields on the associated
bundle with this local form. When $\F_p$ is inhomogeneous, its
tangent module cannot be generated by linear fields alone.

The bundle $P\to L$ also relates two notions of holonomy. Its
component monodromy gives Fernandes's reduced holonomy for the
linearized algebroid (Proposition~\ref{prop:fernandes-holonomy})
and the formal outer holonomy of the maximal tangent module. Around
a closed leaf, Proposition~\ref{prop:formal-maximal} proves local
finite generation of the normal formal completion of the maximal
tangent sheaf, even when the infinitesimal leaves are not closed.
Proposition~\ref{prop:outer-comparison} identifies the image of
formal outer holonomy with the component group of the linear
holonomy isotropy, $\operatorname{im}\Xi\cong H/H^0=\Gamma$.
Consequently, a compact leaf $L$ admits a strong stability model
with compact transverse leaves if and only if $\F_p$ is closed and
$\operatorname{im}\Xi$ is finite (Corollary~\ref{cor:formal-reeb}).
In terms of the formal neighbourhood theory of
Fischer--Laurent-Gengoux~\cite{FischerLaurentGengoux2026}, the
infinitesimal foliation describes the transverse geometry, while
formal outer holonomy detects the component monodromy of the bundle.

\subsection*{Organization of the paper}
Section~\ref{sec:setting} recalls the geometric and function-theoretic
background, the stability notions, and the semi-local groupoid and its
Lie algebroid. Section~\ref{sec:closed} proves the homogeneity and
properness criteria and the associated-bundle model around closed
leaves. It ends by treating properness around a closed saturated
submanifold contained in a stratum. Section~\ref{sec:reeb} gives the compact-leaf criterion,
derives its structural consequences, and applies it to maps and
singular fibers. Global stability follows in Section~\ref{sec:global}.
Section~\ref{sec:modules} compares reduced algebroid holonomy and
formal outer holonomy with component monodromy.

\subsection*{Acknowledgements}

This paper was written for the SIGMA special issue
\emph{Interactions of Poisson Geometry, Lie Theory and Symmetry}
celebrating the 60th birthday of Rui Loja Fernandes. M.~de Melo warmly thanks Rui for his support and many stimulating
conversations, and especially for introducing him to Molino's theory.
Rui's guidance has had a lasting influence on his mathematical career. I.S. is grateful to Rui Loja Fernandes, his doctoral advisor,
collaborator, and friend, for his guidance and for the many
conversations that have shaped the way he thinks about geometry.
It is a pleasure to contribute to this special issue in Rui's honor.

\section{Background}
\label{sec:setting}

In this section we recall the basic definitions and results concerning the theory of singular Riemannian foliations and Lie groupoids which are used in the paper. We also review the construction of the semi-local model \cite{AIMS2022} and discuss the stability problem in foliation theory and how it is related to the linearization of Lie groupoids.

\subsection{Singular Riemannian foliations}

Let $(M,\metric)$ be a Riemannian manifold and $\F$ a partition of $M$ into connected immersed submanifolds, called leaves. Write $L_x$ for the leaf through $x$ and $\mathfrak X(\F)$ for the space of all smooth vector fields on $M$ tangent to the leaves.

\begin{definition}[Singular Riemannian foliation]
The partition $\F$ is a \emph{singular Riemannian foliation} if
\begin{enumerate}[label=(\roman*)]
\item for every $x\in M$ and $v\in T_xL_x$, there is an $X\in\mathfrak X(\F)$ with $X_x=v$,
\item every geodesic perpendicular to one leaf remains perpendicular to each leaf it meets.
\end{enumerate}
\end{definition}

The first condition defines a smooth singular foliation. The second says that the metric is adapted to it. The foliation is \emph{regular} if all its leaves have the same dimension. The connected components of the union of leaves of any given dimension are the \emph{strata} of $\F$. Each stratum $\Sigma_\alpha$ is an embedded submanifold, and $\F_\alpha=\F|_{\Sigma_\alpha}$ is a regular Riemannian foliation for the induced metric \cite[Section~6.2 and Proposition~6.3]{Molino1988}.

A map between foliated manifolds is \emph{foliated} if it sends each leaf into a leaf. A \emph{foliated diffeomorphism} is a diffeomorphism that sends each leaf onto a leaf. A subset is \emph{saturated} if it is a union of leaves.

\begin{example}[Isometric actions]
The connected components of the orbits of an isometric Lie group action form a singular Riemannian foliation \cite[Proposition~3.78]{AlexandrinoBettiol2015}. Such foliations are called \emph{homogeneous}.
\end{example}

For a regular foliation, following plaques along a leafwise path from $p$ to $q$ gives a germ of a diffeomorphism between local transversals at $p$ and $q$. This germ depends only on the homotopy class of the path relative to its endpoints. The germs of loops based at $p$ form the \emph{holonomy group} $\Hol_p(\F)$. The covering of $L_p$ corresponding to the kernel of $\pi_1(L_p,p)\to\Hol_p(\F)$ is its \emph{holonomy cover} \cite[Section~2.1]{MoerdijkMrcun2003}.

\begin{example}[Suspension]\label{ex:suspension}
Let $\widetilde B$ be the universal cover of a connected manifold $B$, and let $\varphi\colon\pi_1(B)\to\operatorname{Diff}(T)$ be a homomorphism. The diagonal action
\[
\gamma\cdot(\widetilde b,u)=(\gamma\widetilde b,\varphi(\gamma)u)
\]
defines the suspension $M_\varphi=(\widetilde B\times T)/\pi_1(B)$. The images of $\widetilde B\times\{u\}$ are the leaves of a regular foliation \cite[Chapter~1]{MoerdijkMrcun2003}. For \mbox{$\Gamma_u=\{\gamma:\varphi(\gamma)u=u\}$}, the leaf through $[\widetilde b,u]$ is diffeomorphic to $\widetilde B/\Gamma_u$, and its holonomy consists of the germs at $u$ of elements of $\varphi(\Gamma_u)$. If $T$ is Riemannian and $\varphi$ acts by isometries, the product of the lifted metric of $B$ and the metric of $T$ descends to a metric adapted to the suspension foliation.
\end{example}

Let $P$ be a sufficiently small plaque and $\operatorname{Tub}_r(P)$ a normal tube on which the normal exponential map is a diffeomorphism. The leaves of the restricted foliation are the connected components of the intersections with $\operatorname{Tub}_r(P)$. Throughout, components of intersections of a leaf with a tube or slice are taken in the intrinsic topology of the leaf.

\begin{lemma}[Homothetic transformation lemma {\cite[Lemma~6.2]{Molino1988}}]\label{lem:homothetic}
For $0<\lambda\leq1$, the map
\[
h_\lambda\colon\operatorname{Tub}_r(P)\longrightarrow
\operatorname{Tub}_{\lambda r}(P),\qquad
\exp_q(v)\longmapsto\exp_q(\lambda v),
\]
is a foliated diffeomorphism.
\end{lemma}

Write $B_\varepsilon(V)$ for the open ball of radius $\varepsilon$ about the origin in a Euclidean space $V$, and $B_\varepsilon(E)$ for the corresponding disk bundle of a Euclidean vector bundle $E$. For $p\in P$ and $V=\nu_pL_p$, pull back the connected components of the intersections of leaves with the slice $\exp_p(B_\varepsilon(V))$ to $B_\varepsilon(V)$. Lemma~\ref{lem:homothetic} extends this partition by positive homotheties to all of $V$. The resulting \emph{infinitesimal foliation} $\F_p$ is a singular Riemannian foliation for the Euclidean metric. The origin is a leaf and each other leaf lies in a sphere centered at the origin \cite[Chapter~6, Proposition~6.5]{Molino1988}.

\subsubsection*{Molino's sheaf}
For a complete singular Riemannian foliation, the regular locus is open and dense in each connected component. On each stratum $\Sigma_\alpha$, the Molino sheaf $\mathscr C_\alpha$ of the regular foliation $\F_\alpha$ is a locally constant sheaf of Lie algebras of germs of transverse Killing fields. A transverse Killing field is locally represented on transversals by Killing fields compatible with holonomy.

\begin{theorem}[{\cite[Theorem~9.8]{AlexandrinoCaramello2022}}]
\label{thm:molino-structure}
Let $(M,\F)$ be a complete singular Riemannian foliation. The partition $\overline\F$ by leaf closures is a singular Riemannian foliation for the same metric. The Molino sheaf of the regular locus extends to a locally constant sheaf $\mathscr C_\F$ on $M$. Restriction to each stratum induces a surjection
\[
\mathscr C_\F|_{\Sigma_\alpha}\longrightarrow\mathscr C_\alpha.
\]
\end{theorem}

The statement about leaf closures is the Molino--Alexandrino--Radeschi (MAR) theorem \cite{AlexandrinoRadeschi2017}. An alternative proof is given in \cite{deMeloStruchiner2026}. Each leaf closure lies in the stratum of the leaf. For $x\in\Sigma_\alpha$, the transverse evaluations of germs in $(\mathscr C_\alpha)_x$ span $T_x\overline{L_x}/T_xL_x$ \cite[Chapter~6, Lemma~6.4]{Molino1988}, \cite[Section~3.1]{Caramello2024}.

\subsubsection*{Basic functions and polynomial presentations}
A function is \emph{basic} if it is constant on each leaf. For an infinitesimal foliation $(V,\F_0)$ with closed leaves, the algebra $\R[V]^{\F_0}$ of basic polynomials is finitely generated and separates leaves \cite[Theorem~1.1 and Section~4]{LytchakRadeschi2018}. If $\rho_1,\ldots,\rho_k$ generate this algebra, every smooth basic function on $V$ is of the form $a(\rho_1,\ldots,\rho_k)$ for some $a\in C^\infty(\R^k)$ \cite[Theorem~B]{MendesRadeschi2019}. The distance between leaves defines a metric on $V/\F_0$ \cite[Section~2.1]{MendesRadeschi2020}.

We also use the partition obtained by grouping together the leaves in each orbit of a finite group acting on $V/\F_0$ by maps induced by orthogonal foliated transformations. For this partition, finite generation, separation and smooth factorization still hold \cite[Lemmas~10--11]{MendesRadeschi2019}.

\subsection{Lie groupoids}\label{subsection: Lie groupoids}

For a Lie groupoid $\G\rra M$, write $s,t\colon\G\to M$ for source and target and $gh$ for the product of arrows with $s(g)=t(h)$. Manifolds are second countable. The base and source fibers are Hausdorff, while the arrow manifold may be non-Hausdorff. The isotropy group at $x$ is $\G_x=s^{-1}(x)\cap t^{-1}(x)$ and its orbit is $\mathcal O_x=t(s^{-1}(x))$. The connected components of the orbits form the \emph{orbit foliation} \cite[Corollary~3.8 and Theorem~8.2]{Meinrenken2017}. They are the orbits of the open source-connected subgroupoid $\G^\circ$ obtained by taking the component of $1_x$ in each $s^{-1}(x)$ \cite[Proposition~1.29]{CrainicFernandes2011}. For $S\subset M$, write $\G|_S=s^{-1}(S)\cap t^{-1}(S)$ for the restriction whenever it is a Lie subgroupoid.

\begin{example}[Action groupoid]\label{ex:action-groupoid}
An action of $\G$ on $N$ with moment map $\mu\colon N\to M$ defines the action groupoid $\G\ltimes N\rra N$. Its arrows are $(g,z)\colon z\to gz$ with $s(g)=\mu(z)$, and $(g,hz)(h,z)=(gh,z)$. In particular, a Lie group action gives $K\ltimes N\rra N$. A \emph{representation} is an action on a vector bundle by fiberwise linear isomorphisms. It is \emph{orthogonal} if these preserve the fiber metrics.
\end{example}

\begin{example}[Holonomy groupoid]\label{ex:holonomy-groupoid}
For a regular foliation $\F$, the arrows of $\Hol(\F)\rra M$ are leafwise paths modulo equality of endpoints and holonomy germs. Multiplication is concatenation, with $\beta\alpha$ traversing $\alpha$ first. Its orbits are the leaves, its isotropy at $p$ is $\Hol_p(\F)$, and its source fiber at $p$ is the holonomy cover of $L_p$ \cite[Section~5.2]{MoerdijkMrcun2003}.
\end{example}

\begin{example}[Submersion groupoid]\label{ex:submersion-groupoid}\label{ex:simple}
For a submersion $f\colon M\to N$, the submersion groupoid $M\times_N M\rra M$ has one arrow from $x$ to $y$ whenever $f(x)=f(y)$. Its orbits are the fibers, and their connected components are the leaves of the simple foliation $\F_f$. This foliation has trivial holonomy. The endpoint map identifies $\Hol(\F_f)$ with $(M\times_N M)^\circ$, and hence with $M\times_N M$ when the fibers are connected \cite[Example~1.14]{CrainicStruchiner2013}.
\end{example}

\begin{example}[Gauge groupoid]\label{ex:gauge-groupoid}
For a principal $H$-bundle $\pi\colon P\to L$, the arrow $[q',q]$ of the gauge groupoid $(P\times P)/H\rra L$ represents the equivariant map $qh\mapsto q'h$ from the fiber over $\pi(q)$ to the fiber over $\pi(q')$. Every transitive Lie groupoid $\G\rra L$ has this form. Fix $p\in L$ and put $H=\G_p$ and $P=s^{-1}(p)$. Then $t\colon P\to L$ is a principal $H$-bundle, and $[q',q]\mapsto q'q^{-1}$ identifies its gauge groupoid with $\G$ \cite[Theorem~3.10]{Meinrenken2017}.
\end{example}

\subsubsection*{Properness and linearization}
\begin{definition}[{\cite[Definition~1.1]{CrainicStruchiner2013}}]
A Lie groupoid $\G\rra M$ is \emph{proper} if its arrow manifold is Hausdorff and $(s,t)\colon\G\to M\times M$ is proper. It is \emph{$s$-proper} if its arrow manifold is Hausdorff and $s$ is proper.

We say that $\G$ is \emph{proper at $x$} if every sequence $g_n$ with $(s(g_n),t(g_n))\to(x,x)$ has a convergent subsequence. It is \emph{$s$-proper at $x$} if the same holds whenever $s(g_n)\to x$. We say that $\G$ is \emph{proper at $L$} if it is proper at every $x\in L$. For the linear holonomy groupoid $\Gl$, these points are the zero vectors $0_x$, $x\in L$.
\end{definition}

A transitive groupoid as in Example~\ref{ex:gauge-groupoid} is proper exactly when $H$ is compact, and $s$-proper exactly when $P$ is compact \cite[Example~1.5]{CrainicStruchiner2013}.

For an embedded orbit $\mathcal O$, the restriction $\G|_{\mathcal O}$ acts on $\nu\mathcal O=TM|_{\mathcal O}/T\mathcal O$ by the normal representation. An arrow $g\colon x\to y$ sends $[v]$ to $[dt_g(w)]$, where $ds_g(w)=v$. Its action groupoid is the \emph{normal groupoid} $\nu_{\mathcal O}(\G)=\G|_{\mathcal O}\ltimes\nu\mathcal O$. A \emph{linearization} is an isomorphism between restrictions of $\G$ and $\nu_{\mathcal O}(\G)$ to neighbourhoods of $\mathcal O$ and the zero section, respectively, which restricts to the identity on $\G|_{\mathcal O}$ \cite[Section~1.2 and Definition~1.10]{CrainicStruchiner2013}.

\subsubsection*{Stability}
\label{subsec:stability}

\begin{definition}[{\cite[Definition~4.2]{CrainicStruchiner2013}}]
\label{def:stability}
A leaf is \emph{stable} if it admits arbitrarily small saturated open neighbourhoods.
\end{definition}

The classical local Reeb stability theorem gives a model for such neighbourhoods around a compact leaf $L$ of a regular foliation with finite holonomy. Fix $p\in L$, put $H=\Hol_p(\F)$ and $V=\nu_pL$, and let $\widetilde L\to L$ be the holonomy cover. Then $L$ has arbitrarily small saturated neighbourhoods foliated diffeomorphic to
\[
\widetilde L\times_H V,
\]
where $H$ acts on $V$ by linear holonomy. Indeed, the finite holonomy action on a transverse disk can be linearized, and a ball for an $H$-invariant inner product is radially $H$-equivariantly diffeomorphic to $V$. The leaves are the images of $\widetilde L\times\{v\}$ and are finite covers of $L$ \cite{Reeb1952}, \cite[Section~2.3]{MoerdijkMrcun2003}.

\begin{theorem}[{\cite[Theorem~1 and Corollary~2]{CrainicStruchiner2013}}]
\label{thm:groupoid-linearization}
Let $\G\rra M$ be Hausdorff and let $\mathcal O$ be the orbit through $x$. If $\G$ is proper at $x$, it admits a linearization
\[
\G|_W\cong\nu_{\mathcal O}(\G)|_{W'}
\]
for neighbourhoods $W$ of $\mathcal O$ and $W'$ of the zero section. If $\G$ is $s$-proper at $x$, there are arbitrarily small invariant neighbourhoods $W$ of $\mathcal O$ for which $\G|_W\cong\nu_{\mathcal O}(\G)$.
\end{theorem}

When $\Hol(\F)$ is Hausdorff, the Reeb model follows by applying this theorem to the holonomy groupoid. Its source fibers are connected, and the source fiber at $p$ is the compact holonomy cover of $L$. Thus $\Hol(\F)$ is $s$-proper at $p$, and its normal representation is linear holonomy \cite[Remark~1.2 and Example~1.15]{CrainicStruchiner2013}.

For a singular foliation we use the following associated-bundle model. Let $P\to L$ be a principal $H$-bundle with connected total space, where $H$ is a compact Lie group acting linearly on $V$. Let $\F_0$ be a smooth singular foliation on $V$, invariant under positive homotheties, with $\{0\}$ as a leaf. Assume that $H$ preserves $\F_0$ and that $H^0$ preserves each leaf. For an $H$-invariant, relatively compact, saturated open neighbourhood $\Omega\subset V$ of $0$, the images of $P\times L'$ for leaves $L'$ of $\F_0|_\Omega$ define the foliation on $P\times_H\Omega$.

\begin{definition}[Strong stability]\label{def:reeb-stability}
A compact leaf $L$ is \emph{strongly stable} if it admits a \emph{local associated-bundle model}, that is, a foliated diffeomorphism
\begin{equation}\label{eq:reeb-stability-model}
(U,\F|_U)\cong P\times_H(\Omega,\F_0|_\Omega)
\end{equation}
from a saturated open neighbourhood of $L$ to a model as above, carrying $L$ to the zero section. The model has \emph{compact transverse leaves} if every leaf of $\F_0|_\Omega$ is compact.
\end{definition}

In Theorem~\ref{thm:closed}, the transverse foliation is $\F_p$ on $V=\nu_pL$, the action is orthogonal, and we take $\Omega=B_\varepsilon(V)$.

\subsection{Semi-local model}\label{sec: semi-local}

Let $B$ be a closed saturated submanifold contained in a stratum of $(M,\F)$, and put $E=\nu B$. We work on a saturated tube $U$ carrying the semi-local model of \cite{AIMS2022}, identified by the normal exponential map with a neighbourhood of the zero section in $E$. When we use a constant-radius tube, we assume one is available. For a closed leaf in a complete manifold this follows from \cite[Proposition~16 and Theorem~A]{MendesRadeschi2019}. Closedness of leaves in a tube is understood relative to the tube.

Write $\pi\colon E\to B$ for the projection and $\F_B=\F|_B$. For $b\in B$, pull back the connected components of the intersections of leaves with the normal slice at $b$ and extend by positive homotheties to $E_b$. This gives the \emph{reduced infinitesimal foliation} $\F_b$. When $B=L$ is a leaf, it is the infinitesimal foliation on $\nu_bL$.

In bundle coordinates, $h_\lambda(v)=\lambda v$. The linearization along $B$ of a vector field $X$ tangent to $\F$ is
\[
X^\ell(v)=\lim_{\lambda\to0^+}(dh_\lambda)_v^{-1}X(h_\lambda(v)).
\]
It is smooth and tangent to $\F$, and its flow acts by linear isometries between fibers. The flows of these vector fields generate the \emph{linearized foliation} $\Fl$. In particular, $\Fl\subset\F$ and $\Fl|_B=\F_B$ \cite[Sections~2.1.1--2.1.2]{AIMS2022}.

Extend $T\F_B$ to a distribution tangent to $\F$ on $U$ and linearize it along $B$. The resulting horizontal distribution $\mathcal T$ is tangent to $\Fl$, and $d\pi$ maps $\mathcal T_v$ isomorphically onto $T_{\pi(v)}\F_B$. It defines a metric partial connection $\nabla^\tau$ along $\F_B$ \cite[Proposition~2.2]{AIMS2022}. Parallel transport $P_\alpha$ along a leafwise path $\alpha$ in $B$ preserves the fiber foliations.

For $b\in B$, put
\[
\kk_b=\{A\in\so(E_b): Av\in T_v L_v^E\text{ for every }v\in E_b\},
\]
where $L_v^E$ is the leaf of $\F_b$ through $v$. The connected Lie subgroup $K_b^0\subset O(E_b)$ with Lie algebra $\kk_b$ is the largest connected group of linear isometries preserving every leaf of $\F_b$ \cite[Example~2.1]{AIMS2022}. Thus $\F_b$ is homogeneous exactly when its leaves are the $K_b^0$-orbits.

Let $\F^E$ denote the fiberwise foliation with leaves $L_v^E$. In the semi-local model, $\F$ is generated by $\F^E$ and horizontal transport \cite[Theorem~1.4 and Proposition~3.1]{AIMS2022}. More explicitly, for $v\in U$ and $b=\pi(v)$, its leaf and its linearized leaf are
\begin{equation}\label{eq:semilocal}
\begin{aligned}
L_v      &=\{P_\alpha(w):w\in L_v^E,\ \alpha\in\mathcal P_b\},\\
L_v^\ell &=\{P_\alpha(kv):k\in K_b^0,\ \alpha\in\mathcal P_b\},
\end{aligned}
\end{equation}
where $\mathcal P_b$ consists of the piecewise smooth paths in the leaf of $\F_B$ through $b$ that start at $b$. In particular, $T_vL_v=\mathcal T_v\oplus T_vL_v^E$.

To realize $\Fl$ as an orbit foliation, let $Q=\Fr(E)$ be the orthonormal frame bundle and $n=\operatorname{rank}E$. Each linearized flow $\varphi_t$ lifts to $Q$ by $\widetilde\varphi_t(\xi)=\varphi_t|_{E_b}\circ\xi$ for a frame $\xi\colon\R^n\to E_b$. These lifts commute with the right $O(n)$-action. The foliation $\widetilde\F$ they generate is regular. At $\xi$ its tangent space is the sum of the horizontal lift of $T_b\F_B$ and the tangent space to $K_b^0\xi$. The action on frames is free, and $\dim K_b^0$ is locally constant \cite[Corollary~2.5 and proof of Theorem~4.1]{AIMS2022}.

The right $O(n)$-action lifts to a free proper action on $\Hol(\widetilde\F)$ by groupoid automorphisms. The quotient
\[
\G=\Hol(\widetilde\F)/O(n)\rra B
\]
acts orthogonally on $E$. Indeed, for $g\colon b\to b'$ and $\xi\in Q_b$ there is a unique representative $\widetilde g$ with source $\xi$. The action is given by $g\cdot\xi(u)=t(\widetilde g)(u)$ for $u\in\R^n$, independently of the frame. Fixing $\xi$ identifies the source fiber of $\G$ at $b$ with the source fiber of $\Hol(\widetilde\F)$ at $\xi$, so it is connected \cite[Proposition~2.14 and proof of Theorem~4.1]{AIMS2022}.

The \emph{linear holonomy groupoid} is the action groupoid $\Gl=\G\ltimes E\rra E$. Its orbits are the leaves of $\Fl$ \cite[Theorem~4.1]{AIMS2022}. We also write $\Gl$ for its restriction to the tube. For a regular Riemannian foliation around a closed leaf $L$, this construction gives $\Hol(\F)|_L\ltimes\nu L$ \cite[Example~4.5]{AIMS2022}.

When $B=L$ is a closed leaf, $\G$ is transitive. Fix $p\in L$ and put
\begin{equation}\label{eq:notation}
V=E_p,\qquad P=s^{-1}(p),\qquad H=\G_p,\qquad \Gamma=H/H^0.
\end{equation}
Write $\rho\colon H\to O(V)$ for the isotropy representation. The map $[q,v]\mapsto q\cdot v$ identifies $P\times_H V$ with $E$ \cite[Section~1.2]{CrainicStruchiner2013}. The isotropy of $\Gl$ at $0_p$ is $H$.

\subsubsection*{The Lie algebroid over a leaf}
For $B=L$ a closed leaf, the connection $\nabla^\tau$ and the bundle $\mathfrak k\to L$ determine a transitive Lie subalgebroid $A=A(\nabla^\tau,\mathfrak k)$ of the gauge algebroid $\mathfrak{gl}(E)$, fitting into the exact sequence
\[
0\longrightarrow\mathfrak k\longrightarrow A\longrightarrow TL\longrightarrow0.
\]
The connection splits this sequence as vector bundles. In the splitting $A\cong TL\oplus\mathfrak k$, the inclusion into the derivations of $E$ gives the flat $A$-connection
\[
D_{(X,\phi)}s=\nabla^\tau_Xs+\phi(s)
\]
\cite[Definition~6.4 and Proposition~6.5]{AIMS2022}.

In fact $A=\operatorname{Lie}(\G)$. After choosing a frame, the source fiber of $\G$ is the holonomy cover of its lifted leaf. The tangent space of that leaf consists of horizontal lifts of $TL$ and vertical $\mathfrak k$-directions. Differentiating the action on frames identifies $\operatorname{Lie}(\G)$ with the corresponding derivations of $E$ \cite[Proposition~6.7 and Remark~6.8]{AIMS2022}. Hence
\[
\mathcal A^\ell=A\ltimes E=\operatorname{Lie}(\Gl)
\]
has orbit foliation $\Fl$. Its transverse algebroid at $0_p$, on the normal fiber $V=E_p$, is
\begin{equation}\label{eq:transverse-action-algebroid}
\mathcal T_p=\mathfrak k_p\ltimes V.
\end{equation}
With the commutator bracket on $\mathfrak k_p$, a section $u\colon V\to\mathfrak k_p$ has anchor $X_u(v)=-u(v)v$. The orbits of $\mathcal T_p$ are the $K_p^0$-orbits.

\section{Properness around a closed leaf}
\label{sec:closed}
We first ask when $\F$ agrees near the closed saturated submanifold
$B$ with its linearization along $B$. In that case, the orbits of
$\Gl$ in the tube are precisely the leaves of $\F|_U$.

\begin{proposition}\label{prop:linearizable}
In the semi-local model around $B$, $\F|_U=\Fl|_U$ if and
only if $\F_b$ is homogeneous for every $b\in B$. If $B$ is a leaf closure, homogeneity at one point suffices.
\end{proposition}

\begin{proof}
Fix $b\in B$ and write $E_b=\nu_bB=(T_bB)^\perp$. In normal
coordinates, the leaves of $\F_b$ in $U\cap E_b$ are the connected
components of the intersections of the leaves of $\F$ with this slice.
For $\Fl$, the corresponding components are the $K_b^0$-orbits.
If $\F|_U=\Fl|_U$, these partitions agree on the slice and hence,
by positive homotheties, on all of $E_b$. Thus $\F_b$ is homogeneous.

Conversely, suppose that the leaves of $\F_b$ are the orbits of a
connected subgroup $G_b\subset O(E_b)$. Then $G_b\subset K_b^0$,
and for every $v\in E_b$ we have
\[
G_bv\subset K_b^0v\subset L_v^E=G_bv.
\]
If this holds for every $b\in B$, the two descriptions in
\eqref{eq:semilocal} agree.

Finally, suppose that $B$ is a leaf closure. Linearized leafwise
fields, together with the foliated linear extensions of
$\F_B$-basic fields supplied by
\cite[Proposition~2.4]{AIMS2022}, span the base directions locally.
Composing their flows gives local foliated linear trivializations
of the reduced infinitesimal foliations. These flows are fiberwise
linear and preserve $\F$, hence also the connected components of
its intersections with the fibers. The maps, initially defined
on small fiber balls, extend by homotheties.
Since $B$ is connected, composition along a path gives a foliated
linear isomorphism
\[
P_{bc}\colon (E_b,\F_b)\longrightarrow(E_c,\F_c)
\qquad(b,c\in B).
\]
If $\F_b$ is given by a connected subgroup $G_b\subset O(E_b)$,
then $\F_c$ is given by $G_c=P_{bc}G_bP_{bc}^{-1}$.
Every element of $G_c$ preserves each leaf of $\F_c$, and these
leaves lie in spheres centred at the origin. Hence $G_c$ preserves
the norm and lies in $O(E_c)$. Thus homogeneity at one point implies
homogeneity at every point of $B$, and the preceding argument applies.
\end{proof}

\begin{corollary}\label{cor:diffeomorphic}
Let $L$ be a closed leaf and $p\in L$. There is a foliated
diffeomorphism between $\F$ and $\Fl$ near $L$, preserving $L$, if and only if $\F_p$ is homogeneous.
\end{corollary}

\begin{proof}
If $\F_p$ is homogeneous, Proposition~\ref{prop:linearizable}
gives $\F=\Fl$ on a tube, so the identity is such a diffeomorphism.
Conversely, let $\phi$ carry $\F$ to $\Fl$ and put $q=\phi(p)\in L$. In a local product chart for $\Fl$ one may use the product of a plaque metric and the Euclidean metric on the fiber. By \cite[Proposition~2.4]{AlexandrinoRadeschi2017}, the normal differential $D^\nu_p\phi$ carries $\F_p$ onto the infinitesimal foliation of $\Fl$ at $q$, whose leaves are the $K_q^0$-orbits. Conjugating $K_q^0$ by $(D^\nu_p\phi)^{-1}$ gives a connected group with orbits the leaves of $\F_p$. Since these lie in spheres about zero, every element preserves the norm on $\nu_pL$. Thus $\F_p$ is homogeneous.
\end{proof}

\begin{example}[The octonionic Hopf foliation]\label{ex:octonionic-linearization}
The Clifford foliation associated with the system $C_{8,1}$ on
$\R^{16}$ in Radeschi's construction is the octonionic Hopf
fibration $S^{15}\to S^8$
\cite[Theorem~A and Corollary~2.9]{Radeschi2014}.
In octonionic coordinates, it is the restriction of the quadratic map
\[
h\colon \Oct^2\longrightarrow\R\oplus\Oct,
\qquad h(x,y)=\bigl(|x|^2-|y|^2,\,2x\bar y\bigr).
\]
Since $|h(x,y)|=|x|^2+|y|^2$, its fibers on $\R^{16}$ are
$\{0\}$ and the radial copies of the Hopf fibers, each diffeomorphic
to $S^7$. They define a singular Riemannian foliation
$\F_{\Oct}$ for the Euclidean metric, the conical extension of
Radeschi's Clifford foliation. It is invariant under positive
homotheties and has the origin as its only point leaf, so it is
an infinitesimal foliation. Its restriction to $S^{15}$ is regular
but inhomogeneous \cite[Theorem~A(2)]{Radeschi2014}, and
$\F_{\Oct}$ is inhomogeneous as well.

To illustrate Proposition~\ref{prop:linearizable}, in this example every linear vector field tangent to $\F_{\Oct}$ is zero
\cite[Lemma~4.19]{NahariStrobl2024}. Thus $K_0^0=\{e\}$,
and linearization along $\{0\}$ gives
\[
\F_{\Oct}^{\ell}=\bigl\{\{v\}:v\in\R^{16}\bigr\}
\ne\F_{\Oct}.
\]
In particular, the two foliations cannot be foliated
diffeomorphic near $0$, since their nonzero leaves have
different dimensions.
\end{example}

\begin{remark}
Notice that the linear holonomy groupoid of Example \ref{ex:octonionic-linearization} is the unit groupoid of $\R^{16}$, which is proper. Hence the homogeneity assumption in Proposition~\ref{prop:linearizable} cannot be replaced by properness of $\Gl$, even when all leaves of $\F$ are compact. This is directly related to our next discussion and Reeb stability theorem (Theorem \ref{thm:reeb}). 
\end{remark}

We now relate closedness of the leaves to properness of the linear holonomy groupoid through its isotropy. Let $L$ be a closed leaf and use the notation
\eqref{eq:notation}. With $V=\nu_pL$, write
\begin{equation}\label{notation along leaf}
q\colon P\times B_\varepsilon(V)\longrightarrow P\times_H B_\varepsilon(V),
\qquad J=\overline{\rho(H)}\subset O(V).
\end{equation}
Here $q$ is the quotient map, and we identify its target with $U$.

\begin{lemma}\label{lem:structure}
The representation $\rho\colon H\to O(V)$ is faithful and
$\rho(H^0)=K_p^0$. Moreover, $P$ is connected and
\[
L_v\cap B_\varepsilon(V)=\rho(H)L_v^E
\qquad(v\in B_\varepsilon(V)).
\]
\end{lemma}

\begin{proof}
Fix a frame $\xi\in Q_p$.
Each $h\in H$ has a unique representative $\widetilde h\in\Hol(\widetilde\F)$ with source $\xi$ by the construction in \cite[Theorem~4.1]{AIMS2022}. Therefore, an element of $\ker\rho$ is represented by a holonomy
arrow of $\widetilde\F$ from $\xi$ to itself.
Its holonomy is represented by a composition of lifted
linearized flows, by construction (see section \ref{sec: semi-local}), which gives that the resulting linear map on $E_p$
fixes the frame $\xi$ and is therefore the identity.
Moreover, it lifts to the identity on $\Fr(E)_p$.
Since $\widetilde\F$ projects onto the single leaf $L$,
this fiber contains a transversal to $\widetilde\F$, because the horizontal decomposition projects isomorphically to the distribution $T\F|_L$.
The holonomy germ is consequently trivial. This proves faithfulness.

For a one-parameter subgroup $h_t$ of $H$, the curve
$\rho(h_t)v$ is the target of the path $(h_t,v)$ in a source
fiber of $\Gl$, since it is an action groupoid. It is therefore smooth in the intrinsic manifold
structure of a leaf of $\Fl\subset\F$. Since the curve is vertical,
its velocity is tangent to $\F_p$, by the tangent space description
following \eqref{eq:semilocal}.
Thus $\rho(H^0)\subset K_p^0$.
Conversely, the $K_p^0$-orbits on $\Fr(E)_p$ lie in the
leaves of $\widetilde\F$
\cite[Section~4, equation~(5)]{AIMS2022}.
Paths in these orbits give paths of arrows starting
at the identity, so $K_p^0\subset\rho(H^0)$. This gives the second assertion.

Finally, every arrow of $\Gl$ is represented by compositions of
linearized flows, which preserve the leaves of $\F$.
Conversely, the parallel transports in
\eqref{eq:semilocal} are induced by arrows of $\G$.
Thus $L_v$ is the $\Gl$-saturation of $L_v^E$.
Arrows starting and ending over $p$ act through $H$, giving
$L_v\cap B_\varepsilon(V)=\rho(H)L_v^E$. The connected source fiber of $\Hol(\widetilde\F)$ at $\xi$ maps onto $P$ under the quotient, so $P$ is connected.
\end{proof}
\begin{remark}
Since $P$ is connected, the principal $\Gamma$-cover $P/H^0\to L$,
with $\Gamma=H/H^0$, is connected. Its monodromy therefore gives
a surjection
\begin{equation}\label{eq:components}
\partial\colon\pi_1(L,p)\longrightarrow\Gamma.
\end{equation}
In particular, $\Gamma$ is finite whenever $\pi_1(L,p)$ is finite.
\end{remark}

\begin{lemma}\label{lem:closed-K}
If every leaf of $\F_p$ is closed, then $K_p^0$ is
closed in $O(V)$ and hence compact.
\end{lemma}

\begin{proof}
The closure of $K_p^0$ is connected.
Since limits of isometries preserving a closed leaf
still preserve it, this closure preserves every leaf
of $\F_p$. Maximality of $K_p^0$ gives the conclusion.
\end{proof}

\begin{theorem}\label{thm:closed}
Let $L$ be a closed leaf and $p\in L$. Then
$\Fl|_U$ is closed if and only if $\Gl|_U$ is proper.
Moreover, $\F|_U$ is closed if and only if $\Gl|_U$ is proper
and $\F_p$ is closed. Consequently, if $\F|_U$ is closed, then
\[
(U,\F|_U)\cong P\times_H(B_\varepsilon(V),\F_p),
\]
where $H$ is compact and $P\to L$ is a principal $H$-bundle.
If, in addition, $L$ is compact, it is strongly stable in the sense of
Definition~\ref{def:reeb-stability} and hence stable.
\end{theorem}

\begin{proof}
The representation of $\G$ on $E$ identifies $U$ with
$P\times_H B_\varepsilon(V)$ by $q(\eta,v)\mapsto\eta\cdot v$.
For any subset $C\subset B_\varepsilon(V)$, the relation
$q(\eta h,v)=q(\eta,\rho(h)v)$ gives
\[
q^{-1}\bigl(q(P\times C)\bigr)=P\times\rho(H)C.
\]
In particular, if $C$ is $H$-invariant, then $q(P\times C)$ is
closed in $U$ if and only if $C$ is closed in $B_\varepsilon(V)$,
since $q$ is a quotient map.

Identify $B_\varepsilon(V)$ with the fiber of $U$ over $p$.
By Lemma~\ref{lem:structure}, the original leaf through $v$ meets
this fiber in $S_v=\rho(H)L_v^E$, the union of the $H$-translates
of its infinitesimal leaf. The linearized leaf meets the same fiber
in the single orbit $N_v=\rho(H)v$. The original leaf is the
$\Gl$-saturation of $L_v^E$, and the linearized leaf is the
$\Gl$-orbit of $v$. Thus
\[
\begin{aligned}
L_v      &=q(P\times L_v^E)=q(P\times S_v),\\
L_v^\ell &=q(P\times\{v\})=q(P\times N_v).
\end{aligned}
\]
The sets $q(P\times L_v^E)$ are connected because $P$ and $L_v^E$
are connected. They are the leaves of the associated-bundle foliation,
so the identification above gives a foliated diffeomorphism
\[
(U,\F|_U)\cong P\times_H(B_\varepsilon(V),\F_p).
\]
It follows from the quotient criterion that $L_v$ is closed exactly
when $S_v$ is closed, and $L_v^\ell$ is closed exactly when $N_v$
is closed. Both $S_v$ and $N_v$ lie in the sphere of radius $|v|$,
so their closedness is equivalent to compactness.

We then relate properness to the isotropy group.
The restriction of $\Gl$ to the zero section is $\G$.
Since $\G$ is transitive, properness at $p$ is equivalent
to compactness of $H$. In that case $\G$ is proper.
The action groupoid of a proper groupoid is proper.
In particular, properness at $0_p$, properness on $U$,
and compactness of $H$ are equivalent.

By Lemma~\ref{lem:structure}, $\rho$ restricts to a bijective Lie
group homomorphism from $H^0$ onto $K_p^0$, with its intrinsic Lie
group structure. Such a homomorphism is an isomorphism
\cite[Proposition~1.1.8]{MonclairNotes}. Hence $H^0$ is compact if and
only if $K_p^0$ is compact, that is, closed in $O(V)$.
Since the components of a Lie group are open cosets of its identity
component, it is compact precisely when that component is compact
and the component group is finite.
We therefore have
\begin{equation}\label{eq:proper-isotropy}
\Gl|_U\text{ is proper}
\quad\text{if and only if}\quad
K_p^0\text{ is closed in }O(V)\text{ and }|\Gamma|<\infty.
\end{equation}

We now prove the assertion about the linearized leaves.
If $H$ is compact, every $N_v$ is compact, so $\Fl|_U$ is closed.
Conversely, suppose every $N_v$ is compact.
A locally closed orbit of an action of a Lie group is an embedded submanifold
\cite[Theorem~4.4.2]{MonclairNotes}.
The $H^0$-orbits, which are the $K_p^0$-orbits, are
the connected components of $N_v$.
They are therefore compact, and there are finitely many.

The compact group $J$ preserves every $N_v$, and $J^0$ preserves each
connected component of $N_v$.
Thus $J^0$ preserves every $K_p^0$-orbit and hence
every leaf of $\F_p$. It follows that
\[
J^0\subset K_p^0=\rho(H^0)\subset\rho(H)\subset J.
\]
The subgroup $\rho(H)$ contains the open subgroup
$J^0$, so it is open and closed in $J$.
It is dense by definition of $J$, hence $\rho(H)=J$.
Faithfulness and \cite[Proposition~1.1.8]{MonclairNotes} identify
$H$ with the compact group $J$. This proves the first assertion.

It remains to characterize closedness of the original leaves.
If $H$ and every infinitesimal leaf are compact,
then every $S_v$ is compact, so the leaves of $\F|_U$
are closed.

Conversely, suppose every $S_v$ is compact.
Since $\rho(H^0)=K_p^0$ preserves each infinitesimal leaf,
$S_v$ is a countable union of infinitesimal leaves,
indexed by $\Gamma$, which is countable since the components of the second-countable Lie group $H$ are open.
All these leaves have the same dimension.

The closure of $L_v^E$ is contained in $S_v$ and is a connected
submanifold, foliated by leaves of $\F_p$.
Indeed, the Euclidean space $V$ is complete, so the leaf-closure theorem applies to $\F_p$ by \cite[Chapter~6]{Molino1988} and
\cite{AlexandrinoRadeschi2017}.
If this closure had dimension greater than $L_v^E$, the countably many leaves covering
$\overline{L_v^E}$ would have measure zero in it, which is
impossible. Here the leaf inclusions are smooth,
and the leaves are second countable.
Consequently the closure and $L_v^E$ have the same dimension.
The leaves in $\overline{L_v^E}$ are then open in
$\overline{L_v^E}$, so connectedness gives
$\overline{L_v^E}=L_v^E$.
Thus every leaf of $\F_p$ is closed.

The infinitesimal leaf space is now a metric space
\cite[Section~2.1]{MendesRadeschi2020}.
The image of $S_v$ in it is compact and countable,
and $\rho(H)$ acts transitively on this image.
A nonempty countable compact metric space has an
isolated point by Baire's theorem.
Transitivity makes every point isolated, and
compactness implies that there are finitely many.
Hence $S_v$ is a finite union of infinitesimal leaves.

The group $J$ preserves $S_v$.
Since its finitely many infinitesimal leaves are
disjoint compact subsets, $J^0$ preserves each of them.
This holds for every $v$, so $J^0\subset K_p^0$.
The preceding open-subgroup argument again shows
that $H$ is compact.

Finally, suppose that $\F|_U$ is closed. The second assertion shows
that $\Gl|_U$ is proper, hence that $H$ is compact. Together with the
foliated diffeomorphism $(U,\F|_U)\cong P\times_H(B_\varepsilon(V),\F_p)$
obtained at the beginning of the proof, this gives the principal
bundle model, with compact structure group $H$. The remaining
conditions of the associated-bundle model follow from
Lemma~\ref{lem:structure} and closedness of $\F_p$, whose
leaves lie in distance spheres. If $L$ is compact, the model
gives strong stability by Definition~\ref{def:reeb-stability}.
\end{proof}

\begin{corollary}\label{cor:finite-components}
Suppose that $\F_p$ is closed. Let $H$ be the isotropy group of $\Gl$
at $0_p$. Then $\F|_U$ is closed if and only if the component group
$\Gamma=H/H^0$ is finite. In particular, if $\pi_1(L,p)$ is finite,
then $\F|_U$ is closed. If $L$ is compact and $\Gamma$ is finite,
then $L$ is strongly stable.
\end{corollary}

\begin{proof}
Apply Lemma~\ref{lem:closed-K}, Theorem~\ref{thm:closed}, and the
surjection $\partial\colon\pi_1(L,p)\to\Gamma$ in \eqref{eq:components}.
\end{proof}

\begin{corollary}\label{cor:homogeneous-proper}
If $\F_p$ is homogeneous, then $\Gl|_U$ is proper if and only
if $\F|_U$ is closed.
\end{corollary}

\begin{proof}
Combine Proposition~\ref{prop:linearizable} with
Theorem~\ref{thm:closed}.
\end{proof}

\subsection*{Closed saturated submanifolds}

For a single closed leaf $L$, the groupoid $\G\rra L$ is transitive,
so its properness is determined by the isotropy at one point.
Now let $B$ be a closed saturated submanifold contained in a stratum.
Recall from Section~\ref{sec: semi-local} that $\G\rra B$ acts on
$E=\nu B$ and that $\Gl=\G\ltimes E$. The orbits of $\G$ are the
leaves of $\F|_B$. When $B$ contains several leaves, compact isotropy
alone no longer guarantees properness.

\begin{proposition}\label{prop:proper-B}
For the linear holonomy groupoid around $B$, $\Gl$ is proper
if and only if $\G$ is proper.
Either condition implies that $\F_B$ is closed, and holds
whenever $\Hol(\widetilde\F)$ is proper.
\end{proposition}

\begin{proof}
The restriction of $\Gl$ to the zero section is $\G$,
and the action groupoid of a proper groupoid is proper.
This proves the equivalence.
The orbits of $\G$ are the leaves of $\F_B$, and
proper groupoids have closed orbits, so $\F_B$ is closed.

Finally, the frame bundle projection $\pi_Q\colon Q\to B$ is proper because its fiber
$O(n)$ is compact.
For a compact subset $C\subset B\times B$, its inverse
image in $Q\times Q$ is compact.
Properness of $\Hol(\widetilde\F)$ makes the inverse
image of this set under its source-target map compact.
Its image in the quotient
$\G=\Hol(\widetilde\F)/O(n)$ is exactly the inverse image
of $C$ under the source-target map of $\G$.
The quotient arrow manifold is Hausdorff, since $O(n)$ is compact
and acts on the Hausdorff arrow manifold of $\Hol(\widetilde\F)$.
Thus $\G$ is proper, and so is $\Gl$.
\end{proof}

In particular, if $B$ is a leaf closure and $\Gl$ is proper,
then its dense leaf is closed in $B$ and hence equals $B$.
A useful sufficient condition for
properness of $\Hol(\widetilde\F)$ is that
$\widetilde\F$ consist of the connected fibers of a
submersion $f\colon Q\to N$ onto a Hausdorff manifold.
Then
\[
\Hol(\widetilde\F)\cong Q\times_NQ
\]
is a closed subgroupoid of the pair groupoid and is
proper.

The distinction from the single-leaf case is already
visible for a regular foliation.
For the product of an irrational linear foliation
of $T^2$ with the foliation of $\R^k$ by points,
take $B=T^2\times\{0\}$.
The groupoid over $B$ has trivial isotropy but dense
orbits, so it is not proper.

\begin{corollary}\label{cor:proper-B-model}
If every reduced infinitesimal foliation $\F_b$
is homogeneous and $\G\rra B$ is proper, then
$\F|_U=\Fl|_U$ is the orbit foliation of the proper
groupoid $\Gl|_U$.
\end{corollary}

\begin{proof}
Apply Proposition~\ref{prop:linearizable} and
Proposition~\ref{prop:proper-B}.
\end{proof}

\section{Reeb stability}
\label{sec:reeb}

For compact leaves, an associated-bundle model with compact
transverse leaves also has compact nearby leaves. Combined with
Theorem~\ref{thm:closed}, this gives the following characterization
and a criterion for strong stability.

\begin{theorem}[Reeb stability]\label{thm:reeb}
Let $L$ be a compact leaf and $p\in L$. After shrinking the saturated
tube $U$ around $L$, the following are equivalent:
\begin{enumerate}[label=(\roman*)]
\item $L$ admits a strong stability model with compact transverse leaves.
\item Every leaf of $\F|_U$ is compact.
\item $\Gl$ is proper at $L$ and $\F_p$ is closed.
\end{enumerate}
\end{theorem}

\begin{proof}
Suppose (i) holds. The principal bundle in the stability model has
compact base and compact structure group, so its total space is
compact. Each model leaf is the image of its product with a compact
transverse leaf. A sufficiently small saturated tube about $L$ lies
in this model, which proves (ii).

Each leaf in the tube lies at a fixed distance from $L$. Since $L$
is compact, every sphere bundle of radius less than the tube radius
is compact. Thus closedness in the tube is equivalent to compactness
for its leaves. By Theorem~\ref{thm:closed}, closedness is equivalent
to properness of $\Gl|_U$ together with closedness of $\F_p$.
Properness on $U$ and properness at $L$ are both equivalent to
compactness of the isotropy group $H$, as shown in the proof of that
theorem. This proves the equivalence of (ii) and (iii).

Under (iii), the associated-bundle model of Theorem~\ref{thm:closed}
has compact structure group. The leaves of $\F_p$ are closed and
lie in Euclidean spheres, so they are compact. This gives (i).
\end{proof}

Under these equivalent conditions, $L$ is stable. Moreover, the
principal bundle $P\to L$ in Theorem~\ref{thm:closed} has compact
total space, so its gauge groupoid $\G=(P\times P)/H$ is Hausdorff
and $s$-proper. The same holds for $\Gl|_U$, since $U$ is invariant
and the source map of the action groupoid is the pullback of the
source map of $\G$.

\begin{example}[Classical Reeb stability]\label{rem:classical}
Let $\F$ be a regular Riemannian foliation and $L$ a compact leaf.
Then $\F_p$ is the foliation of $V=\nu_pL$ by points. It is closed
and homogeneous, so $\F|_U=\Fl|_U$ by
Proposition~\ref{prop:linearizable}, and the linear holonomy
groupoid is the normal groupoid of the holonomy groupoid along $L$,
\[
\Gl=\nu_L(\Hol(\F))=\Hol(\F)|_L\ltimes\nu L
\]
\cite[Example~4.5]{AIMS2022}. Its isotropy $H$ at $0_p$ is the
holonomy group $\Hol_p(\F)$, with its discrete Lie group topology,
acting faithfully on $V$ by linear holonomy
(Lemma~\ref{lem:structure}), and $P$ is the holonomy cover
$\widetilde L\to L$. Thus $\Gl$ is proper at $L$ exactly when
$\Hol_p(\F)$ is finite, and in this case Theorem~\ref{thm:reeb} gives
\[
(U,\F|_U)\cong\widetilde L\times_{\Hol_p(\F)} B_\varepsilon(V),
\]
whose leaves are the images of the sets $\widetilde L\times\{v\}$.
This is the suspension of the linear holonomy homomorphism
$\pi_1(L,p)\to\Hol_p(\F)\subset O(V)$ (Example~\ref{ex:suspension}),
that is, the classical local Reeb stability model
\cite[Section~2.3]{MoerdijkMrcun2003}. When the holonomy groupoid is
Hausdorff, it is also the model given by its invariant linearization
\cite[Example~1.15]{CrainicStruchiner2013}. The leaf through
$q(\eta,v)$ is $\widetilde L/\Hol_p(\F)_v$, where $\Hol_p(\F)_v$ is
the stabilizer of $v$, and it covers $L$ with
$[\Hol_p(\F):\Hol_p(\F)_v]$ sheets. Hence, in classical Reeb
stability, the leaves are in fact all compact in the tube on which
the linearization is carried out, and the equivalence in
Theorem~\ref{thm:reeb} reduces to the statement that this happens
exactly when the holonomy of $L$ is finite.
\end{example}

The associated-bundle description also identifies $\Gl|_U$ up to
Morita equivalence with $H\ltimes B_\varepsilon(V)$: the equivalence
is implemented by $P\times B_\varepsilon(V)$, with the commuting
actions induced by the gauge groupoid of $P$ and by $H$.

Even for a regular Riemannian foliation, a compact leaf can be stable
while the nearby leaves are noncompact.

\begin{example}[A stable leaf with infinite holonomy]\label{ex:irrational}
Let $R$ be a rotation of $\R^2$ through an irrational multiple of
$2\pi$, and consider the suspension of
\[
\varphi\colon \pi_1(S^1)\cong\Z\longrightarrow O(2),
\qquad \varphi(n)=R^n,
\]
as in Example~\ref{ex:suspension}. Thus
\[
M_\varphi=(\R\times\R^2)/\Z,\qquad
n\cdot(t,v)=(t+n,R^nv),
\]
with leaves the images of $\R\times\{v\}$.
Since $0$ is fixed by $\varphi$, it determines a compact leaf
$L\cong S^1$, with holonomy $\langle R\rangle\cong\Z$.
For $v\ne0$ the stabilizer is trivial, so the corresponding leaf is
diffeomorphic to $\R$ and dense in the torus
$(\R\times S^1_{|v|})/\Z$.
The suspension is a regular Riemannian foliation and equals its
linearization along $L$, but $\Gl$ is not proper because its
holonomy is infinite. The images of $\R\times B_\varepsilon(\R^2)$ form
a basis of saturated neighbourhoods of $L$, so $L$ is stable.
It is not strongly stable: in a regular foliation the transverse
foliation in an associated-bundle model is by points, and $H^0$
acts trivially. A compact structure group would therefore give
finite holonomy, contrary to $\Hol(L)\cong\Z$.
\end{example}

We now ask when the leaves near $L$ form the orbit foliation of a Hausdorff
Lie groupoid which is proper at $L$. For a Lie groupoid $\G'$,
write $(\G')^\circ$ for its source-connected subgroupoid.

\begin{corollary}\label{cor:proper-model}
Let $L$ be a compact leaf and $p\in L$. Near $L$, the foliation $\F$ admits a Hausdorff
Lie groupoid model which is proper at $L$ if and only if $\F_p$
is homogeneous and $L$ admits a strong stability model with compact
transverse leaves.
\end{corollary}

\begin{proof}
If $\F_p$ is homogeneous and $L$ admits such a strong stability
model, Theorem~\ref{thm:reeb} gives compactness of nearby leaves, and
Proposition~\ref{prop:linearizable} gives $\F=\Fl$.
Theorem~\ref{thm:reeb} then shows that $\Gl$ is a Hausdorff
$s$-proper model. Equivalently, one may assume homogeneity
of $\F_p$ and properness of $\Gl$: in that case
$K_p^0=\rho(H^0)$ is compact, so its orbits are compact,
and Theorem~\ref{thm:reeb} applies.

Suppose now that $\G'$ is Hausdorff, proper at $L$,
and has orbit foliation $\F$.
The source-connected subgroupoid $(\G')^\circ$ has orbit $L$ through $p$, and its isotropy
$((\G')^\circ)_p$ is an open, hence closed, subgroup of the
compact group $\G'_p$. Thus $((\G')^\circ)_p$ is compact.
The source fiber of $(\G')^\circ$ over $p$ is a principal
$((\G')^\circ)_p$-bundle over $L$, hence compact.

Since the source fiber of $(\G')^\circ$ over $p$ is compact and the source fibers of $(\G')^\circ$ are connected,
we obtain that $(\G')^\circ$ is $s$-proper
at $p$ (this follows from Ehresmann theorem, for instance, check
\cite[Remark 1.2]{CrainicStruchiner2013}).
Invariant linearization (Theorem~\ref{thm:groupoid-linearization})
then identifies $(\G')^\circ$ near $L$ with its normal model
on an invariant neighbourhood of the zero section.

In a local bundle chart of the normal model, choose an isotropy-invariant inner product on the transverse fiber. Its infinitesimal leaves are the orbits of the identity component of $((\G')^\circ)_p$. The normal differential of the foliated diffeomorphism identifies these orbits with the leaves of $\F_p$ by \cite[Proposition~2.4]{AlexandrinoRadeschi2017}. The conjugated action is orthogonal by the sphere argument in Corollary~\ref{cor:diffeomorphic}.
Its leaves in a sufficiently small invariant tube
are compact, since both $L$ and $((\G')^\circ)_p$ are compact.
The same holds for $\F$, and Theorem~\ref{thm:reeb} gives a strong
stability model with compact transverse leaves, as required.
\end{proof}

For Corollaries~\ref{cor:finite-cover} and~\ref{cor:leaf-bundle}, let $L$ be a closed leaf and suppose
that $\F|_U$ is closed. Theorem~\ref{thm:closed} supplies the model
\[
(U,\F|_U)\cong P\times_H(B_\varepsilon(V),\F_p),\qquad H\text{ compact},\qquad
\Gamma=H/H^0\text{ finite}.
\]
The bundle $P\to L$ is connected and $H^0$ preserves each
leaf of $\F_p$. 
For compact $L$, these assumptions are equivalent to the conditions of
Theorem~\ref{thm:reeb}.

Consider the quotient and the pullback tube
\[
\widehat L=P/H^0,\qquad
\widehat U=\widehat L\times_L U\cong P\times_{H^0} B_\varepsilon(V),
\]
and let $\widehat\F$ be the foliation by connected components
of inverse images of leaves under $\widehat U\to U$.

\begin{corollary}\label{cor:finite-cover}
The map $\widehat L\to L$ is a connected finite regular
covering with deck group $\Gamma$, and
\[
(\widehat U,\widehat\F)
 =P\times_{H^0}(B_\varepsilon(V),\F_p),
\qquad
\widehat U/\widehat\F\cong B_\varepsilon(V)/\F_p.
\]
\end{corollary}

\begin{proof}
The principal $H$-bundle
$P\to L$ makes $P/H^0\to L$ a principal $\Gamma$-bundle,
and it is connected because $P$ is connected. Its pullback of
$U$ is $P\times_{H^0} B_\varepsilon(V)$, and the covering map is
$[\eta,v]_{H^0}\mapsto[\eta,v]_H$.

The connected leaves of the lifted foliation are
$P\times_{H^0}L_v^E$. Indeed, these are connected, have
the lifted tangent distribution, and together partition
$\widehat U$. Thus $[\eta,v]\mapsto[L_v^E]$ induces the
stated homeomorphism of quotient spaces. The residual deck
action gives
\[
U/\F\cong(\widehat U/\widehat\F)/\Gamma.
\]
\end{proof}

For $v\in B_\varepsilon(V)$, write $L_v$ for the corresponding
leaf in $U$, set $\pi_v=\pi|_{L_v}$, and put
\[
\Gamma_{[v]}=\{\gamma\in\Gamma:\gamma[L_v^E]=[L_v^E]\}.
\]

We use $\rho(\gamma)L_v^E$ to mean the image under any
representative of $\gamma\in\Gamma$ in $H$. This is well defined
because $H^0$ preserves each infinitesimal leaf. Thus
$S_v=\bigcup_{\gamma\in\Gamma}\rho(\gamma)L_v^E$.

\begin{corollary}\label{cor:leaf-bundle}
The projection $\pi_v\colon L_v\to L$ is a proper submersion with fiber
\[
S_v=\bigcup_{\gamma\in\Gamma}\rho(\gamma)L_v^E.
\]
With compatible basepoints,
\[
[\pi_1(L):\operatorname{im}(\pi_v)_*]
 =[\Gamma:\Gamma_{[v]}].
\]
\end{corollary}

\begin{proof}
The distinct leaves in $S_v$ are disjoint compact
embedded submanifolds, all of the same dimension. They are
exactly the connected components of $S_v$, indexed by
$\Gamma/\Gamma_{[v]}$.
The description
\[
L_v=P\times_H S_v\longrightarrow L
\]
gives local trivializations of $\pi_v$. Its fiber $S_v$ is compact, so
$\pi_v$ is proper.
The boundary map for $P\to L$ gives a surjection
$\pi_1(L)\to\Gamma$. The monodromy action of $\pi_1(L)$ on
$\pi_0(S_v)$ is the induced transitive action on
$\Gamma/\Gamma_{[v]}$. The homotopy exact sequence of the
fiber bundle identifies $\operatorname{im}(\pi_v)_*$ with
the stabilizer of the chosen fiber component. Its index is
therefore the number of components.
\end{proof}

Thus the finite coverings of nearby leaves in regular Reeb
stability become fiber bundles with finitely many connected
infinitesimal leaves as fiber. The fiber-bundle description
also appears in the slice approach of
\cite[proof of Proposition~22]{MendesRadeschi2019}. The formula
above expresses the number of components through the isotropy
of $\Gl$. When $\F_p$ is homogeneous, one may equivalently
write $L_v=P/H_v$, where $H_v$ stabilizes $v$.

\begin{corollary}[Polynomial presentation and smooth factorization]
\label{cor:polynomial}
Let $L$ be a closed leaf and let $U_\varepsilon$ be a tube carrying the
semi-local model. The following are equivalent:
\begin{enumerate}[label=(\roman*)]
\item $\F|_{U_\varepsilon}$ is the partition into connected fiber
components of a continuous map to a Hausdorff space,
\item $\F|_{U_\varepsilon}$ is closed,
\item $\Gl|_{U_\varepsilon}$ is proper and $\F_p$ is closed,
\item and there is an $H$-invariant polynomial map $\tau\colon V\to\R^r$ such
that $\widehat\tau([\eta,v])=\tau(v)$ has exactly the leaves as fibers.
\end{enumerate}
The map $\tau$ can be chosen from generators of
$\R[V]^{\F_p,H}$, the algebra of $H$-invariant $\F_p$-basic
polynomials on $V$. For this choice and $0<\delta<\varepsilon$, every smooth basic $f\colon U_\varepsilon\to\R^m$ satisfies
\[
f|_{U_\delta}=\sigma\circ\widehat\tau|_{U_\delta}
\quad\text{for some }\sigma\in C^\infty(\R^r,\R^m).
\]
Moreover $U_\varepsilon/\F\cong(B_\varepsilon(V)/\F_p)/\Gamma$.
\end{corollary}

\begin{proof}
Fiber components of a continuous Hausdorff-valued map are closed,
so (i) implies (ii). Theorem~\ref{thm:closed} gives (ii) if and only if (iii).
Under (iii), properness gives compactness of $H$, and closedness of
$\F_p$ gives compactness of its leaves. The group $H^0$ preserves each infinitesimal leaf,
and $\Gamma$ is finite. Lemma~\ref{lem:structure} identifies the
intersection of a leaf with $B_\varepsilon(V)$ with $H L_v^E$.
Apply \cite[Lemma~10]{MendesRadeschi2019} to this disconnected
infinitesimal partition to obtain $\tau$. Its invariance makes
$\widehat\tau$ well defined, and its fibers are exactly
$P\times_H(H L_v^E)$. This proves (iv), which implies (i).

A basic $f$ restricts to an $H$-invariant $\F_p$-basic function $g$
on $B_\varepsilon(V)$, with $f([\eta,v])=g(v)$. Multiply each component
of $g$ by a radial cutoff supported in $B_\varepsilon(V)$ and
equal to one on $B_\delta(V)$, then extend by zero. Theorem~B and
Lemma~11 of \cite{MendesRadeschi2019} factor this extension through
the chosen polynomial generators. The displayed equality follows.
The quotient description follows from the same leaf intersections.
\end{proof}

This corollary places the polynomial presentation within the properness
criterion of Theorem~\ref{thm:closed}. Lytchak--Radeschi obtain
algebraicity from basic polynomials \cite{LytchakRadeschi2018}, and
Mendes--Radeschi combine the slice theorem with invariant-function
theory \cite{MendesRadeschi2019}. Here the finite transverse symmetry
group is identified as the component group of the linear isotropy.
Their slice structure group consists of all orthogonal transformations
preserving each disconnected slice leaf, and it contains $\rho(H)$.
In the closed setting their identity components agree, but the full
groups may differ: for the octonionic foliation of
Example~\ref{ex:octonionic-linearization}, $H$ is trivial whereas
$-I$ preserves each Hopf fiber.

As an application, we study stability of maps whose connected
fiber components form singular Riemannian foliations. Stability
and strong stability of a component refer to
Definitions~\ref{def:stability} and~\ref{def:reeb-stability} for
this foliation. Let $f\colon M\to N$ be such a continuous map, with $M$ complete
and $N$ Hausdorff. Every fiber component $L$ is closed, and a
saturated tube around it consists of closed leaves.
Theorem~\ref{thm:closed} supplies an associated-bundle model.
If $L$ is compact, Theorem~\ref{thm:reeb} gives strong stability
with compact nearby leaves, and Corollary~\ref{cor:proper-model}
provides the criterion for a proper Lie groupoid model.

\begin{corollary}[Stability of singular maps]
\label{cor:fibers}
Let $f\colon M\to N$ be a continuous map from a complete Riemannian manifold to a Hausdorff space. Suppose its connected fiber components form a singular Riemannian foliation $\F_f$.
Every compact fiber component is strongly stable and hence stable.
Near such a component $L$, the foliation $\F_f$ admits a Hausdorff Lie groupoid model which is proper at $L$ if and only if its infinitesimal foliation is
homogeneous.
\end{corollary}

The nearby leaf bundles of Corollary~\ref{cor:leaf-bundle}
describe the fiber components in these neighbourhoods.

\begin{example}[A singular circle fiber]\label{ex:hopf-mapping-torus}
Let $h\colon\Oct^2\to\R\oplus\Oct$ be the quadratic Hopf map of
Example~\ref{ex:octonionic-linearization}. Since $h(-v)=h(v)$,
it induces a smooth map
\[
f\colon M=(\R\times\Oct^2)/\Z\longrightarrow\R\oplus\Oct,
\qquad f([t,v])=h(v),
\]
where $n\cdot(t,v)=(t+n,(-1)^n v)$. The product metric descends
to a complete metric on $M$. Locally the fiber foliation is the
product of an interval with the octonionic Hopf foliation, so it
is a singular Riemannian foliation.

The zero fiber is a circle $L$. Each nonzero fiber is the mapping
torus of the antipodal map on a Hopf fiber $S^7$. In particular,
the fibers are connected, and their dimensions jump from one to eight.
Corollary~\ref{cor:fibers} gives strong stability of $L$.
At each $p\in L$ the infinitesimal foliation is the octonionic Hopf
foliation on $\Oct^2$, which is inhomogeneous. Thus no Hausdorff Lie
groupoid which is proper at $L$ has this orbit foliation near $L$.
\end{example}

\section{Global stability}
\label{sec:global}

To pass from the tube of compact leaves given by Theorem~\ref{thm:reeb}
to the whole manifold, we use the Molino sheaf to propagate closedness.

\begin{lemma}\label{lem:propagation}
Let $\F$ be a singular Riemannian foliation on a complete connected
Riemannian manifold $M$. If some nonempty saturated open subset
of $M$ consists of leaves that are closed in $M$, then every leaf
is closed in $M$.
\end{lemma}

\begin{proof}
Let $U$ be such an open set. The regular locus is open and dense,
so its intersection with $U$ is nonempty and open. On this
intersection every leaf is closed. The transverse evaluations of
the Molino sheaf therefore vanish there, and every local section,
being a transverse Killing field, is zero. Since $\mathscr C_{\F}$
is locally constant and $M$ is connected, it vanishes on all of $M$.
Its restriction to each stratum surjects onto the Molino sheaf of
that stratum, whose transverse evaluations describe the leaf closures.
Hence every leaf is closed
by Theorem~\ref{thm:molino-structure} (see also
\cite[Section~9.1 and Theorem~9.8]{AlexandrinoCaramello2022}
and \cite[Section~3.1]{Caramello2024}).
\end{proof}

\begin{theorem}[Global Reeb stability]\label{thm:global}
Let $(M,\F)$ be a singular Riemannian foliation on a compact connected
manifold. The following are equivalent:
\begin{enumerate}[label=(\roman*)]
\item There is a compact leaf $L$ such that $\Gl$ is proper at $L$
and $\F_p$ is closed for some $p\in L$.
\item Some leaf of $\F$ admits a strong stability model with compact
transverse leaves.
\item Every leaf of $\F$ is compact.
\item The leaves are exactly the fibers of a smooth map $M\to\R^N$.
\end{enumerate}
\end{theorem}

\begin{proof}
By Theorem~\ref{thm:reeb}, conditions (i) and (ii) are equivalent,
and either gives a nonempty saturated open tube of compact leaves.
Lemma~\ref{lem:propagation} then implies that every leaf is closed.
Since $M$ is compact, this gives (iii). Conversely, under (iii),
Theorem~\ref{thm:reeb} gives a strong stability model with compact
transverse leaves at every leaf. In particular, (ii) holds and every
leaf is stable.

Assume (iii). Finite generation of smooth basic functions in this compact case also follows from \cite[Theorem~C and the subsequent remark]{MendesRadeschi2019}. To construct a separating map directly, choose finitely many saturated tubes $U_i$ as in
Corollary~\ref{cor:polynomial}, with smaller tubes covering $M$.
On $U_i$ let $\tau_i$ be the resulting smooth map separating leaves.
Choose a basic function $\phi_i$ supported in $U_i$ and positive
on its smaller tube. Such functions are obtained by smooth cutoffs of
the squared normal radius, which is smooth on the tube and constant
on leaves. Extend $\phi_i$ and $\phi_i\tau_i$ by zero outside $U_i$.
Set
\[
\Phi=(\phi_i,\phi_i\tau_i)_i.
\]
This map is smooth and basic. If $\Phi(x)=\Phi(y)$, choose $i$ with
$\phi_i(x)>0$. Equality of the first coordinates gives
$\phi_i(y)=\phi_i(x)>0$, and the remaining coordinates imply
$\tau_i(x)=\tau_i(y)$. Thus $x$ and $y$ are on the same leaf.
This proves (iv). Finally, under (iv), every leaf is a closed subset
of compact $M$, giving (iii).
\end{proof}

\begin{corollary}[Finite fundamental group]\label{cor:global-finite}
Let $(M,\F)$ be a singular Riemannian foliation on a compact connected
manifold. If a compact leaf $L$ has finite fundamental
group and $\F_p$ is closed for some $p\in L$, then all leaves are compact
and are exactly the fibers of a smooth map $M\to\R^N$.
\end{corollary}

\begin{proof}
Corollary~\ref{cor:finite-components} and Theorem~\ref{thm:reeb}
give a saturated neighbourhood of compact leaves.
Apply Theorem~\ref{thm:global}.
\end{proof}

For comparison, Zhukova's global stability theorem implies that, for
a transversely complete orbit-like foliation, a compact regular leaf
with finite fundamental group forces all leaves to be compact
\cite[Theorem~3 and Section~12]{Zhukova2007}.
Corollary~\ref{cor:global-finite} assumes that $M$ is compact and
allows the initial leaf to be singular, provided its infinitesimal
foliation is closed.

\begin{example}[The closedness hypothesis]\label{ex:irrational-sphere}
For $\alpha\in\R\setminus\mathbb Q$, consider the isometric action
of $\R$ on the round sphere
$S^4\subset\C^2\times\R$ given by
\[
t\cdot(z_1,z_2,x)=(e^{it}z_1,e^{i\alpha t}z_2,x).
\]
The poles $(0,0,\pm1)$ are point leaves,
with trivial fundamental group, whereas every orbit with
$z_1z_2\ne0$ is dense in a two-dimensional torus.
At either pole, the infinitesimal foliation is the orbit foliation
of the same irrational linear action on $\C^2$, so it is not closed.
Thus the hypothesis on $\F_p$ in
Corollary~\ref{cor:global-finite} cannot be omitted. The poles are
strongly stable, using the trivial principal bundle over a point and the indicated infinitesimal
foliation as transverse model. Their nonclosed transverse leaves
also show why the compactness of the transverse leaves in Theorem~\ref{thm:global}(ii)
cannot be dropped.
\end{example}

\section{Comparison of holonomy constructions}
\label{sec:modules}
Let $L$ be a closed leaf and fix $p\in L$. Recall that
$\Gamma=H/H^0$ is the component group of the isotropy group $H$
of $\Gl$ at $0_p$. The monodromy of the principal $\Gamma$-cover
$P/H^0\to L$ gives the surjection
$\partial\colon\pi_1(L,p)\to\Gamma$ of \eqref{eq:components}.
We first compare $\partial$ with the reduced holonomy of the action
algebroid. We then compare it with the formal outer holonomy of the
maximal tangent module of $\F$.

\subsection{Reduced Lie algebroid holonomy}
\label{sec:fernandes}
Consider the action algebroid $\mathcal A^\ell$
of Section~\ref{sec: semi-local} and its transverse algebroid
$\mathcal T_p=\mathfrak k_p\ltimes V$ from
\eqref{eq:transverse-action-algebroid}.
Fernandes associates to $L$ the reduced holonomy
\[
\overline{\Hol}^{\mathcal A^\ell}_L\colon
\pi_1(L,p)\longrightarrow\operatorname{Out}(\mathcal T_p),
\]
where outer automorphisms are germs of transverse algebroid
automorphisms modulo flows of time-dependent sections.
The transverse algebroid germ is well defined, and the outer
transport depends only on the homotopy class of the base loop
\cite[Theorem~1.2 and Proposition~3.2]{Fernandes2002}.

Since $\G=(P\times P)/H$, its algebroid is $TP/H$. A principal connection splits $TP/H\to TL$. For a based loop $\gamma$, choose the algebroid path given by this splitting and let $\eta(t)$ be the horizontal lift from $\eta_0=1_p$. Extend $\dot\gamma$ locally in time to a time-dependent base field and lift it horizontally to $P$. Its equivariant flow gives local gauge bisections whose arrow along $\gamma$ is $[\eta(t),\eta_0]$. In Fernandes's construction \cite[Section~3.1]{Fernandes2002}, use this splitting and the fiberwise constant extension of the horizontal section to the action algebroid. Its flow acts on $E$ by the bisection and on transverse isotropy by conjugation. At the endpoint the transverse algebroid map is therefore
\[
\Theta_h(\xi,v)
=\bigl(\rho(h)\xi\rho(h)^{-1},\rho(h)v\bigr),
\qquad \xi\in\mathfrak k_p,
\]
where a lift of the loop ends at its initial point multiplied
by $h\in H$. With the convention
$[\eta h,v]=[\eta,\rho(h)v]$, its component is
$\partial[\gamma]\in\Gamma$, where $\partial$ is the
surjective monodromy of $P/H^0\to L$ from
\eqref{eq:components}.

If $h\in H^0$, then $\Theta_h$ is a product of flows of constant
transverse sections and is inner. Conversely, the anchor of a
section $u\colon V\to\mathfrak k_p$ is $X_u(v)=-u(v)v$, with derivative $D_0X_u=-u(0)\in\mathfrak k_p$ at the origin.
The derivative of an inner flow therefore solves a linear
evolution equation in $\mathfrak k_p$ and lies in $K_p^0$.
If $\Theta_h$ is inner, then
$\rho(h)=D_0\rho(h)\in K_p^0=\rho(H^0)$.
Faithfulness of $\rho$ gives $h\in H^0$.
Thus the kernel of $h\mapsto[\Theta_h]$ is exactly $H^0$,
and the transport description proves the following comparison.

\begin{proposition}\label{prop:fernandes-holonomy}
For a closed leaf $L$, the isotropy representation induces an
injective homomorphism
\[
\iota_{\mathrm{alg}}\colon \Gamma\longrightarrow
\operatorname{Out}(\mathfrak k_p\ltimes V).
\]
With the path convention of Section~\ref{subsection: Lie groupoids},
\begin{equation}\label{eq:fernandes-holonomy}
\overline{\Hol}^{\mathcal A^\ell}_L
=\iota_{\mathrm{alg}}\circ\partial.
\end{equation}
In particular, its image is isomorphic to $\Gamma$.
\end{proposition}

By \eqref{eq:proper-isotropy}, $\Gl|_U$ is consequently proper
exactly when $K_p^0$ is compact and this reduced holonomy image is
finite. For compact $L$ with $\F_p$ closed,
Lemma~\ref{lem:closed-K} and Theorem~\ref{thm:reeb} identify this
finiteness with compactness of nearby leaves, yielding strong
stability in Definition~\ref{def:reeb-stability}.
The transverse point is transversally stable in Fernandes's sense:
Euclidean balls are invariant because $X_u(|v|^2)=0$, so Fernandes's stability theorem applies to
$\mathcal A^\ell$ when $L$ is compact and the reduced holonomy is finite
\cite[Theorem~3.1]{Fernandes2002}.
Its transverse leaves are the $K_p^0$-orbits, and the model for the
original foliation retains the full infinitesimal leaves,
as in Corollary~\ref{cor:leaf-bundle}.
\subsection{The canonical maximal tangent module}
A singular Riemannian foliation is understood here as a smooth
leaf partition. In the literature, the term singular foliation
also refers to an involutive, locally finitely generated
$C^\infty(M)$-submodule of $\mathfrak X_c(M)$
\cite{AndroulidakisSkandalis2009,Hermann1963}.
However, distinct locally finitely generated modules can induce the same leaf
partition \cite[Example~2.9]{NahariStrobl2025} and, moreover, even when a single
vector field generates the partition, the module of all smooth leaf-tangent fields can fail to be locally finitely generated
\cite[Example~2.3]{NahariStrobl2025}.

For infinitesimally closed
singular Riemannian foliations, we prove that the maximal tangent module is locally
finitely generated and inherits the associated-bundle model. We then study the formal tangent module around a closed leaf and use it to compare formal outer holonomy with component monodromy.

For an open subset $U\subset M$, define
\[
\mathscr X_{\F}(U)
=\mathfrak X(\F|_U)
=\{X\in\mathfrak X(U):
X(x)\in T_xL_x\text{ for every }x\in U\}.
\]
Thus $\mathscr X_{\F}$ is the sheaf of tangent fields
introduced in Section~\ref{sec:setting}. Its compactly
supported sections form the module
\[
\mathfrak X_c(\F)
=\mathfrak X(\F)\cap\mathfrak X_c(M)
=\Gamma_c(M,\mathscr X_{\F}).
\]
The sheaf $\mathscr X_{\F}$ depends only on the smooth leaf
partition. It is involutive, and the smoothness axiom for
$\F$ says that its sections span every leaf tangent space.
The additional property discussed here is local finite
generation.
Recall that $\F$ is \emph{infinitesimally closed} if every
infinitesimal foliation $\F_p$ has closed leaves. 
\begin{theorem}
\label{thm:maximal-module}
Let $\F$ be an infinitesimally closed singular Riemannian
foliation. Then $\mathscr X_{\F}$ is locally finitely
generated. Consequently, $\mathfrak X_c(\F)$ is a locally
finitely generated, involutive $C^\infty(M)$-module whose
leaf partition is $\F$. It is maximal among modules
realizing this partition and is natural under foliated
diffeomorphisms.
\end{theorem}

\begin{proof}
Fix $p\in M$. The plaque-local product description and the
homothetic transformation lemma give a foliated coordinate
neighbourhood
\[
W=P_0\times B_\varepsilon(V),
\]
where $P_0\subset\R^d$ is a plaque coordinate neighbourhood
and $B_\varepsilon(V)\subset V=\nu_pL_p$ carries the restriction
of $\F_p$
\cite[Theorem~6.1 and Proposition~6.5]{Molino1988}.
Its local leaves are the connected components of
$P_0\times(L'\cap B_\varepsilon(V))$, where $L'$ is a leaf
of $\F_p$.

Since $\F_p$ is closed, choose basic polynomials
$\rho_1,\ldots,\rho_N$ whose common fibers are its leaves
\cite[Theorem~1.1 and Section~4]{LytchakRadeschi2018}.
Set
\[
F\colon W\longrightarrow\R^N,
\qquad
F(u,v)=(\rho_1(v),\ldots,\rho_N(v)).
\]
Writing $\mathcal C^\infty_W$ for the smooth-function
sheaf, we obtain
\begin{equation}
\label{eq:maximal-kernel}
\mathscr X_{\F}|_W
=
\ker\!\left(
dF\colon (\mathcal C^\infty_W)^{d+\dim V}
\longrightarrow(\mathcal C^\infty_W)^N
\right).
\end{equation}
Indeed, tangent fields annihilate every $\rho_k$.
Conversely, integral curves of a field annihilating $F$
remain in a connected fiber component. These components
are embedded leaf pieces, so the field is leaf-tangent.
This is an equality of smooth solution sheaves, and it does
not require the pointwise kernel of $dF$ to equal the
leaf tangent space at singular points.

Equip $W$ with its standard real-analytic structure.
Since $dF$ has polynomial entries, its analytic solution
sheaf
\[
\mathscr K^\omega
=
\ker\!\left(
dF\colon (\mathcal C^\omega_W)^{d+\dim V}
\longrightarrow(\mathcal C^\omega_W)^N
\right)
\]
is coherent. Near each point, it therefore has finitely
many analytic generators $X_1,\ldots,X_r$.
By \cite[Corollary~VI.1.12]{Malgrange1966}, these fields also generate the smooth solutions. After shrinking the neighbourhood, every smooth field satisfying $dF(X)=0$ can be written as
\[
X=f_1X_1+\cdots+f_rX_r,
\qquad f_j\in C^\infty.
\]
Together with~\eqref{eq:maximal-kernel}, this proves
local finite generation of $\mathscr X_\F$.

Any module realizing $\F$ consists of leaf-tangent
fields and is therefore contained in $\mathfrak X_c(\F)$.
Finally, a foliated diffeomorphism
$\phi\colon (M,\F)\to(N,\F')$ preserves tangency and compact
supports, so
\[
\phi_*\mathfrak X_c(\F)=\mathfrak X_c(\F').
\]
\end{proof}

The proof follows the algebraic strategy used by Nahari--Strobl for the octonionic Hopf foliation \cite[Proposition~4.18 and Appendix~A]{NahariStrobl2024} 
by describing tangent fields by polynomial equations and using Malgrange's flatness theorem to pass to smooth solutions. Here the polynomial equations are supplied by basic polynomials separating the leaves of each closed infinitesimal foliation.

Consider the associated-bundle model
\[
q\colon P\times B\longrightarrow P\times_H B,\qquad
B=B_\varepsilon(V),\quad V=\nu_pL,
\]
where $H$ is compact, and write $\F_{\mathrm{mod}}$ for its foliation.
Upstairs the leaves are the products $P\times L'$, where $L'$ is a
leaf of $\F_p|_B$. Thus a field $Y=(Y_P,Y_B)$ is tangent to the
product foliation precisely when
\[
Y_B(\eta,v)\in T_vL_v^E
\qquad\text{for every }(\eta,v)\in P\times B.
\]
Both components may depend on $(\eta,v)$. Downstairs the maximal
tangent module is the image of the invariant product-tangent fields:
\[
\mathfrak X_c(\F_{\mathrm{mod}})
=q_*\!\left(\mathfrak X_c(P\times\F_p|_B)^H\right).
\]
Here invariance is for the diagonal action
$(\eta,v)h=(\eta h,\rho(h)^{-1}v)$. Two invariant lifts have the
same image precisely when their difference is tangent to the
$H$-orbits. Indeed, $dq$ maps $T_\eta P\oplus T_vL_v^E$ onto the
model leaf tangent space, and $\ker dq$ lies in this sum because
$H^0$ preserves each infinitesimal leaf. A principal connection on
$q$ therefore lifts every model-tangent field to an invariant
product-tangent field. The lift of a compactly supported field has
compact support since $H$ is compact.

\begin{corollary}
\label{cor:module-reeb}
Let $L$ be a closed leaf, $p\in L$, and $U$ a saturated tube carrying
the semi-local model. If $\Gl|_U$ is proper and $\F_p$ is closed,
then $\mathscr X_{\F}|_U$ is locally finitely generated. Under the
foliated diffeomorphism
\[
\Phi\colon (U,\F|_U)\longrightarrow
P\times_H(B_\varepsilon(V),\F_p)
\]
of Theorem~\ref{thm:closed}, its compactly supported sections satisfy
\[
\Phi_*\mathfrak X_c(\F|_U)=\mathfrak X_c(\F_{\mathrm{mod}}).
\]
\end{corollary}

\begin{proof}
In a bundle chart $U_0\times B$, the local leaves are products of
plaques in $U_0$ with infinitesimal leaves. Since $\F_p$ is closed,
the polynomial presentation and kernel argument in the proof of
Theorem~\ref{thm:maximal-module} give local finite generation in
each such chart. The diffeomorphism $\Phi$ preserves leaf tangency
and compact support in both directions, giving the stated
identification.
\end{proof}

If $\F_p$ is inhomogeneous, its maximal tangent module
cannot be generated solely by leaf-tangent linear
fields, since such generators would make its leaves the
$K_p^0$-orbits. The associated-bundle model therefore retains
transverse information beyond the linearized
foliation $\Fl$.
\subsection{Comparison with formal outer holonomy}
Let $L$ be a closed leaf and $p\in L$. Retain the
associated-bundle notation of Lemma~\ref{lem:structure},
with $V=\nu_pL$ and isotropy representation
$\rho\colon H\to O(V)$. Write
\[
K=\rho(H^0)=K_p^0,
\qquad
\kk=\operatorname{Lie}(K).
\]

A \emph{formal singular foliation along $L$} is a
locally finitely generated sheaf of modules of formal
vector fields tangent to $L$, closed under Lie brackets,
whose restrictions span $TL$. In adapted coordinates
$(u,v)$ with $L=\{v=0\}$, the coefficient ring is
$C^\infty(U_0)[[V^*]]$, where coefficients depend smoothly
on the variables along $L$ and formally on the normal
variables. Tangency means that the normal components
vanish at $v=0$.

The \emph{normal formal completion} of $\mathscr X_\F$
along $L$, denoted by $\widehat{\mathscr X}_{\F,L}$,
is the sheaf generated by the normal Taylor expansions
of its smooth sections. Locally, its sections are
finite sums
\[
\sum_j\widehat f_j\operatorname{Tay}_L X_j,
\qquad
\widehat f_j\in C^\infty(U_0)[[V^*]],
\quad X_j\in\mathscr X_\F,
\]
where $\operatorname{Tay}_L X_j$ is obtained by
expanding the coefficients in the normal variable.
No convergence is required. This construction is
independent of the adapted coordinates and discards
fields vanishing to every normal order along $L$.

Similarly, the transverse smooth and formal tangent
modules are
\[
\mathscr M_p=(\mathscr X_{\F_p})_0,
\qquad
\widehat{\mathscr M}_p
=\R[[V^*]]\cdot
\{\operatorname{Tay}_0X:X\in\mathscr M_p\}.
\]
Here $\mathscr M_p$ consists of smooth leaf-tangent
germs at the origin of $V$, and
$\operatorname{Tay}_0X$ denotes their Taylor series.

The following proposition establishes formal finite
generation without assuming properness of $\Gl$ or
closedness of the infinitesimal leaves.

\begin{proposition}
\label{prop:formal-maximal}
\label{lem:maximal-formal}
Let $L$ be a closed leaf of a complete singular
Riemannian foliation $(M,\F)$, with associated-bundle
model $P\times_H(B_\varepsilon(V),\F_p)$.
Then $\widehat{\mathscr X}_{\F,L}$ is a formal singular
foliation along $L$. Its transverse module
$\widehat{\mathscr M}_p$ admits finitely many homogeneous
polynomial tangent generators, and its space of linear
parts is $\kk$.
\end{proposition}

\begin{proof}
Let $\mathscr P_p$ be the module of polynomial vector
fields tangent to $\F_p$. As a submodule of
$\R[V]\otimes V$, it is finitely generated
\cite[Chapter~7]{AtiyahMacdonald1969}.
Homothetic invariance, Lemma~\ref{lem:homothetic}, implies that every homogeneous component of a polynomial tangent field is tangent.
We may therefore choose homogeneous generators
$Q_1,\ldots,Q_r$, of degrees $d_1,\ldots,d_r$.

The same homothetic argument shows that every
homogeneous Taylor term $X^{(d)}$ of a smooth tangent
germ $X$ belongs to $\mathscr P_p$. Indeed, for fixed $v$, the vectors $X(tv)$ for small $t>0$ belong to the fixed subspace $T_vL_v^E\subset V$ under the translation identifications, by homothetic invariance. Taking derivatives at $t=0$ gives the assertion in every degree. By homogeneous
generation,
\[
X^{(d)}=\sum_j f_{j,d-d_j}Q_j,
\]
where $f_{j,d-d_j}$ is homogeneous of degree $d-d_j$
and is zero when $d<d_j$. Collecting these coefficients
degree by degree expresses the Taylor series of $X$
as a formal combination of the $Q_j$. Conversely,
each $Q_j$ is a smooth tangent field. Hence
\[
\widehat{\mathscr M}_p
=\sum_{j=1}^r\R[[V^*]]Q_j.
\]
This module is involutive because brackets of
polynomial tangent fields are tangent.

By the linearization construction in Section~\ref{sec: semi-local}, the linear part of a smooth tangent germ belongs to $\kk$, and every $A\in\kk$ occurs through its fundamental linear field. Since all $Q_j$ vanish at the origin, formal
coefficients introduce no additional linear parts. 

Now choose an associated-bundle chart
$U_0\times B_\varepsilon(V)$. Its tangent fields split
into base directions and vertical fields tangent
to $\F_p$. The preceding degree-by-degree
decomposition preserves smooth dependence on the
base: at each degree, the map from polynomial
coefficients to homogeneous tangent fields is a
surjective map of finite-dimensional vector spaces
and admits a fixed linear right inverse
(compare \cite[Lemma~1.18]{FischerLaurentGengoux2026}).
Consequently,
\[
\widehat{\mathscr X}_{\F,L}
=
\sum_{i=1}^{\dim L}
C^\infty(U_0)[[V^*]]\,\partial_{u_i}
+
\sum_{j=1}^r
C^\infty(U_0)[[V^*]]\,Q_j.
\]
Both inclusions follow from the preceding argument
and the fact that the displayed generators are
smooth tangent fields.
This description proves local finite generation
and involutivity. Naturality of the maximal tangent
sheaf makes the descriptions compatible on overlaps.
Along the zero section, the vertical generators
vanish and the base fields span $TL$. Thus the
completion is a formal singular foliation along
$L$, with transverse module $\widehat{\mathscr M}_p$.
\end{proof}

Let $\operatorname{Sym}(\widehat{\mathscr M}_p)$ be the group
of formal diffeomorphisms fixing zero and preserving this
formal module. Let
$\operatorname{Inner}(\widehat{\mathscr M}_p)$ be the
subgroup generated by its formal flows, and write
\[
\operatorname{Out}(\widehat{\mathscr M}_p)
=
\operatorname{Sym}(\widehat{\mathscr M}_p)
/
\operatorname{Inner}(\widehat{\mathscr M}_p),
\]
as in~\cite[Section~1.4]{FischerLaurentGengoux2026}.

\begin{proposition}
\label{prop:outer-comparison}
Let $L$ be a closed leaf, and use the notation above.
The isotropy representation induces an injective homomorphism
\[
\iota\colon \Gamma=H/H^0
\longrightarrow
\operatorname{Out}(\widehat{\mathscr M}_p).
\]
The formal outer holonomy $\Xi$ of the normal formal completion
of the maximal tangent sheaf along $L$ satisfies
\[
\Xi=\iota\circ\partial,
\qquad
\pi_1(L,p)\xrightarrow{\partial}\Gamma
\xrightarrow{\iota}
\operatorname{Out}(\widehat{\mathscr M}_p),
\]
where $\partial$ is the monodromy of the covering $P/H^0\to L$.
In particular, $\operatorname{im}\Xi\cong\Gamma$, and
\[
\Gl|_U\text{ is proper}
\quad\text{if and only if}\quad
K\text{ is closed in }O(V)
\text{ and }\operatorname{im}\Xi\text{ is finite}.
\]
\end{proposition}

\begin{proof}
Every $\rho(h)$ is a foliated linear isometry of $(V,\F_p)$.
It preserves the maximal tangent module by naturality and
therefore preserves its formal completion.
Thus $H$ acts by formal symmetries of
$\widehat{\mathscr M}_p$.
Every element of $K=\rho(H^0)$ is a product of flows
of fundamental linear fields in $\mathscr M_p$,
so $H^0$ acts by inner symmetries. Hence the action induces a homomorphism $\iota\colon H/H^0\to
\operatorname{Out}(\widehat{\mathscr M}_p)$.
This descent to the component group is the mechanism described in \cite[Example~2.9]{FischerLaurentGengoux2026}.

Conversely, Proposition~\ref{prop:formal-maximal}
and the first-jet argument preceding
Proposition~\ref{prop:fernandes-holonomy} give
\[
D_0\bigl(
\operatorname{Inner}(\widehat{\mathscr M}_p)
\bigr)\subset K.
\]
Thus, if the linear symmetry $\rho(h)$ is inner, then
\[
\rho(h)=D_0\rho(h)\in K=\rho(H^0).
\]
Faithfulness of $\rho$ gives $h\in H^0$.
Hence the kernel of the map from $H$ to outer symmetries
is exactly $H^0$, proving injectivity of $\iota$.
Choose a principal connection on $P$. In a local trivialization of $P\times_HV$, its induced horizontal lifts are base directions plus fundamental linear fields in $\kk$. They are therefore sections of the maximal tangent sheaf. The lift $X\mapsto X^{\mathrm{hor}}$ is $C^\infty(L)$-linear and projects to $X$. Its Taylor expansion thus splits the projection from projectable formal tangent fields to $\mathfrak X(L)$, as required in \cite[Example~1.28]{FischerLaurentGengoux2026}. Formal outer holonomy is independent of this choice \cite[Corollary~1.31, property~$\mathcal P_C$, and Section~2.1]{FischerLaurentGengoux2026}.
A lifted loop has endpoints differing by an element
$h\in H$, and its transverse transport is the formal expansion
of $\rho(h)$.
Its outer class is $\iota(\partial[\gamma])$.
Since $P$ is connected, $\partial$ is surjective by
\eqref{eq:components}, so
$\operatorname{im}\Xi\cong\Gamma$
The properness equivalence follows from
\eqref{eq:proper-isotropy}.
\end{proof}
\begin{corollary}
\label{cor:formal-reeb}
Let $L$ be a compact leaf, $p\in L$, and $U$ a saturated
tube carrying the semi-local model.
Let $\Xi$ be the formal outer holonomy of the normal
formal completion of the maximal tangent sheaf along $L$.
The leaf $L$ admits a strong stability model with compact transverse
leaves if and only if $\F_p$ is closed and
$\operatorname{im}\Xi$ is finite. Equivalently, after shrinking $U$,
every leaf of $\F|_U$ is compact.
\end{corollary}

\begin{proof}
Proposition~\ref{prop:outer-comparison} identifies
$\operatorname{im}\Xi$ with $\Gamma$.
When $\F_p$ is closed, Lemma~\ref{lem:closed-K} and
\eqref{eq:proper-isotropy} identify finiteness of $\Gamma$
with properness of $\Gl|_U$.
The assertion follows from Theorem~\ref{thm:reeb}.
\end{proof}

Propositions~\ref{prop:outer-comparison} and
\ref{prop:fernandes-holonomy} also identify the kernels:
\[
\ker\Xi
=
\ker\partial
=
\ker\overline{\Hol}^{\mathcal A^\ell}_L.
\]
For a comparison of Fischer--Laurent-Gengoux's outer
holonomy with Fernandes's holonomy, see
\cite[Remark~2.5]{FischerLaurentGengoux2026}.
Under the equivalent conditions of Corollary~\ref{cor:formal-reeb}, the connected finite covering $P/H^0\to L$ of Corollary~\ref{cor:finite-cover} is the covering corresponding to the kernel of either formal outer holonomy or reduced algebroid holonomy.

For foliations induced by associated group actions,
the factorization of outer holonomy through the component
group already appears in
\cite[Example~2.9]{FischerLaurentGengoux2026}.
Proposition~\ref{prop:outer-comparison} identifies this
map injectively for the full transverse maximal module,
including inhomogeneous infinitesimal foliations.
The formal classification of
\cite[Theorems~2.16 and~3.1]{FischerLaurentGengoux2026}
combines outer holonomy with principal-bundle data.
The principal $H$-bundle $P\to L$ and its orthogonal action
provide a concrete presentation of the completed maximal
tangent sheaf.

Outer holonomy can retain information lost on the leaf space:
the composite
\[
\pi_1(L,p)\xrightarrow{\partial}\Gamma
\longrightarrow\operatorname{Bij}(V/\F_p)
\]
can have trivial image even when $\Gamma$ is nontrivial. For example, let $\sigma$ be a reflection of $\R^2$ and suspend its action on the foliation by concentric circles: take $(\R\times\R^2)/((t,v)\sim(t+1,\sigma v))$ with the quotient of the product foliation. Along the zero section $L=S^1$, the infinitesimal leaves are circles and the origin. Here $K_p^0=SO(2)$ and a turn around $L$ acts by $\sigma$, so $\rho(H)=O(2)$ and $\Gamma=\Z/2$. The reflection fixes every circle but its derivative has determinant $-1$, so its formal symmetry is outer (compare \cite[Example~1.19]{FischerLaurentGengoux2026}).
Under the common neighbourhood hypotheses, the first
holonomy of Laurent-Gengoux--Ryvkin agrees with
Dazord's leaf-space holonomy
\cite[Remark~2.30]{LaurentGengouxRyvkin2022}.
Their higher holonomies use the universal Lie
$\infty$-algebroid of the transverse module and incorporate
higher homotopy information, beyond the component
monodromy identified here.

\begin{remark}[Formal neighbourhoods of leaf closures]
\label{rem:formal-leafclosure}
The semi-local model suggests extending the formal
neighbourhood classification of
\cite{FischerLaurentGengoux2026}
from a single leaf to a leaf closure $B$.
In this setting, $\F_B$ may have dense nonclosed leaves,
and transport is governed by a generally nontransitive
groupoid over $B$, rather than by a single principal
isotropy bundle $P\to L$.
The maximal tangent sheaf provides a canonical object
whose normal formal completion along $B$ can be studied.
A natural question is whether suitable outer transport
and bundle data classify these formal neighbourhoods.
The semi-local description developed here supplies
a geometric framework for investigating this question.
\end{remark}

\subsection*{Funding and Declarations}

\noindent\textit{Funding.}
I.S. was supported in part by the S\~ao Paulo Research Foundation
(FAPESP), Brasil, under Grant No.~2022/16097-2.
L.R.S. was supported by the S\~ao Paulo Research Foundation
(FAPESP), Brasil, under Grant No.~2025/27715-7.

\noindent\textit{Use of generative AI.}
The authors used ChatGPT (OpenAI) to edit the manuscript and to help
identify possible gaps in the proofs during the final review.

\end{document}